\documentclass[11pt]{article}
\usepackage[left=1in,right=1in,top=1in,bottom=1in]{geometry}
\usepackage{times}
\usepackage{expl3}
\usepackage{cite}
\usepackage[table]{xcolor}
\usepackage{bm}
\usepackage{multirow}
\usepackage{stackengine} 
\usepackage{enumitem}
\usepackage{hhline}
\usepackage{lipsum}
\usepackage{titlesec}
\usepackage{wrapfig}
\usepackage{epsfig}
\usepackage{graphicx}
\usepackage{amsmath}
\usepackage[title]{appendix}
\usepackage{amssymb}
\usepackage{epstopdf}
\usepackage{boldline}
\usepackage{calligra}
\usepackage{url}
\usepackage{mathrsfs}
\usepackage{blindtext}

\newcommand{\define}{\stackrel{\mbox{\tiny def}}{=}}

\newtheorem{definition}{Definition}
\newtheorem{theorem}{Theorem}

\newtheorem{lemma}{Lemma}

\newtheorem{remark}{Remark}

\usepackage{mathtools}
\usepackage{epstopdf}
\usepackage{balance}
\usepackage{thmtools}
\usepackage{thm-restate}
\usepackage{hyperref}
\usepackage{cleveref}

\usepackage[ruled,vlined]{algorithm2e}
\include{pythonlisting}
\newenvironment{roster}
 {\begin{enumerate}[font=\upshape,label=(\alph*)]}
 {\end{enumerate}}
\newcommand{\ostar}{\mathbin{\mathpalette\make@circled\star}}

\makeatletter
\newcommand{\removelatexerror}{\let\@latex@error\@gobble}
\makeatother
\makeatletter
\newcommand*{\rom}[1]{\expandafter\@slowromancap\romannumeral #1@}
\makeatother

\ExplSyntaxOn
\newcommand\latinabbrev[1]{
  \peek_meaning:NTF . {
    #1\@}%
  { \peek_catcode:NTF a {
      #1.\@ }%
    {#1.\@}}}
\ExplSyntaxOff

\titleclass{\subsubsubsection}{straight}[\subsubsection]

\begin{document}
\vspace{1cm}
\title{Random Multiple Operator Integrals}\vspace{1.8cm}
\author{Shih~Yu~Chang 
\thanks{Shih Yu Chang is with the Department of Applied Data Science,
San Jose State University, San Jose, CA, U. S. A. (e-mail: {\tt
shihyu.chang@sjsu.edu}).
           }}

\maketitle

\begin{abstract}
The introduction of Schur multipliers into the context of Double Operator Integrals (DOIs) was proposed by V. V. Peller in 1985. This work extends theorem on Schur multipliers from measurable  functions to their closure space and generalizes the definition of DOIs to Multiple Operator Integrals (MOIs) for integrand functions as Schur multipliersconstructible by taking the limit of projective tensor product and by taking the limit of integral projective tensor product. According to such closure space construction for integrand functions, we demonstrate that any function defined on a compact set of a Euclidean space can be expressed by taking the limit of the projective tensor product of linear functions. We also generalize previous works about random DOIs with respect to finite dimensional operators, tensors, to MOIs with respect to random operators, which are defined from spectral decomposition perspectives. Based on random MOIs definitions and their properties, we derive several tail bounds for norms of higher random operator derivatives, higher random operator difference and Taylor remainder of random operator-valued functions. 
\end{abstract}


\section{Introduction}\label{sec:Introduction} 


Basic concepts about \emph{Double Operator Integrals}(DOIs) were first mentioned in~\cite{daletskii1965integration} by considering integration and differentiation of functions of Hermitian operators and their applications to the theory of perturbations. It was Birman and Solomyak who established later the theory of double operator integrals~\cite{birman1966double,birman1967double,birman2003double}. Let $(\Lambda_1, A_1)$ and $(\Lambda_2, A_2)$ be spaces with spectral measures $A_1$ and $A_2$ on a HIlbert space $\mathfrak{H}$. Given a bounded measurable function $\psi$ and a operator $X$ on Hilber space $\mathfrak{H}$, the DOI is defined as
\begin{eqnarray}\label{eq:DOI def intro}
\int\limits_{\Lambda_1} \int\limits_{\Lambda_2} \psi(\lambda_1,\lambda_2) dA_1(\lambda_1) X dA_2(\lambda_2),
\end{eqnarray}
where $\psi$ is named as an \emph{integrand function}. We use conventional notation $\mathfrak{S}_p$ to represent the $p$-th Schatten ideal. For every $X \in \mathfrak{S}_1$, we say that the integrand function $\psi$ is a \emph{Schur multiplier} of $\mathfrak{S}_1$ associated with the spectral measures $\Lambda_1$ and $\Lambda_2$ if we have
\begin{eqnarray}\label{eq:DOI def S1 intro}
\int\limits_{\Lambda_1} \int\limits_{\Lambda_2} \psi(\lambda_1,\lambda_2) dA_1(\lambda_1) X dA_2(\lambda_2) \in \mathfrak{S}_1.
\end{eqnarray}
Note that the introduction of Schur multipliers into the context of double operator integrals was proposed by ~\cite{peller1985hankel}. This is an extension of the notion of matrix Schur multipliers. The first contribution of this work is to extend the Schur multipliers~\cite{peller2006multiple} from measurable functions to their closure space. We also extend the definition of DOIs to \emph{Multiple Operator Integrals}(MOIs) for integrand as Schur multipliers constructible by the \emph{limit of projective tensor product, $\hat{\otimes}$} and by the \emph{limit of integral projective tensor product, $\acute{\otimes}$}. Therefore, several Lemmas discussed in this work about MOIs basic properties, e.g., continuity, perturbation, and higher order derivative representations by MOIs, will be generalized to functions constructible by taking the \emph{limit of projective tensor product $(\hat{\otimes})$} of measurable functions, or by taking the \emph{limit of integral projective tensor product $(\acute{\otimes})$} of measurable function~\cite{skripka2019multilinear}. According to such closure space construction for integrand functions, we will be able to represent any function $f(\lambda_1, \lambda_2, \cdots, \lambda_m)$ defined on a compact set of a Euclidean space $\mathbb{R}^m$ by the \emph{limit of projective tensor product} of linear functions. 


The linear operators $A_1, A_2$ and $X$ adopted by Eq.~\eqref{eq:DOI def intro} are assumed to be deterministic. The consideration of linear operators under randomness settings is discussed by Skorohod in~\cite{skorohod2001random}. Two senses of linear random operators are defined there: strong random operators and weak random operators. Hackenbroch generalizes the concept of a strong random operator by replacing the index Hilbert space by a fixed dense subspace~\cite{hackenbroch2009point}. In this work, he proves that the operators with a densely defined adjoint have a unique closed extension with a measurable selection, particularly, including symmetric operators. The classes of random operators considered are all restricted by assuming that there is a dense nonrandom subspace, but the results are extensions of the results on self-adjoint extensions obtained by~\cite{skorohod2001random}. Thang and Quy present results on strongly random operators and bounded strongly random operators between separable Banach spaces~\cite{thang2017spectral}. They show that a bounded strongly random operator can be extended to a continuous linear operator from the set of Banach valued random variables to the set of Banach valued random operators equipped with the topology from convergence in probability. Thang and Quy prove several versions of the spectral theorem, including bounded self adjoint, and more generally bounded normal strongly random operators. These results amplify the results obtained by~\cite{skorohod2001random}. Although several existing works about random DOIs for finite dimensional operators, e.g., tensors, have been studied~\cite{chang2022randomDOI, chang2022randomPDOI}, we first attempt to consider random operators with MOIs in this work. However, the randomness of linear operators is equipped with a different randomness structure compared to aforementioned random linear operators. All operators considered in this work are assumed to have spectral decomposition characterized by eigenvalues and projector spaces (unitary operators formed by eigenspaces), then the randomness of a given linear operator is determined by the random variables of eigenvalues and random unitary operators with Haar measure. 




Besides completeness for theorem on Schur multipliers, the extension definition of MOIs, integrand approximation by linear functions, our other contributions include following: derivation properties of MOIs when the integrand functions constructible by the \emph{limit of projective tensor product, $\hat{\otimes}$} and by the \emph{limit of integral projective tensor product, $\acute{\otimes}$}; definition about random MOIs from spectral decomposition perspectives; application MOIs to derive tail bounds for higher random operator derivatives, higher random operator difference and Taylor remainder of random operator-valued functions. 


The remainder of this paper is organized as follows. In Section~\ref{sec:Completenss for Theorem on Schur multipliers}, we discuss the completeness for Theorem on Schur multipliers and define MOIs with the completeness of integrand functions. We represent any function $f(\lambda_1, \lambda_2, \cdots, \lambda_m)$ defined on a compact set of a Euclidean space $\mathbb{R}^m$ by \emph{limit of projective tensor product} of linear functions in Section~\ref{sec:Integrand Approximation By Linear Functions}. We will discuss MOI properties and random MOI for the integrand function constructible by \emph{limit of integral projective tensor product} in Section~\ref{sec:MOIs Properties and Random MOIs}. Finally, we will apply MOIs to derive tail bounds for higher random operator derivatives, higher random operator difference and Taylor remainder of random operator-valued functions in Section~\ref{sec:Random MOI Applications}.

\section{Completenss for Theorem on Schur multipliers}\label{sec:Completenss for Theorem on Schur multipliers}

The purpose of this section is to define MOI with the completeness of integrand funtions.

\subsection{Double Operator Integrals}\label{sec:DOI}


We use $\Gamma(A_1, A_2)$ to represent the space of Schur multipliers of $\mathfrak{S}_1$ associated with the spectral measures $\Lambda_1$ and $\Lambda_2$. If the integrand function $\psi$ in Eq.~\eqref{eq:DOI def intro} belongs to the \emph{projective tensor product} as $L^{\infty}(A_1) \hat{\otimes} L^{\infty}(A_2)$, i.e., the integrand function $\psi$ can be expressed as
\begin{eqnarray}\label{eq:integrand sum form}
\psi(\lambda_1,\lambda_2) = \sum\limits_{n \geq 0} f_{1,n}(\lambda_1) f_{2,n}(\lambda_2),
\end{eqnarray}
where $f_{1,n} \in L^{\infty}(A_1)$,  $f_{2,n} \in L^{\infty}(A_2)$ and $\sum\limits_{n \geq 0} \left\Vert f_{1,n} (\lambda_1) \right\Vert_{L^{\infty}}\left\Vert f_{2,n} (\lambda_2) \right\Vert_{L^{\infty}} < \infty$; we have $\psi \in \Gamma(A_1, A_2)$. Then, for such integrand function $\psi$, we have
\begin{eqnarray}
\int\limits_{\Lambda_1} \int\limits_{\Lambda_2} \psi(\lambda_1,\lambda_2) dA_1(\lambda_1)  dA_2(\lambda_2) \in \mathfrak{S}_1 = \sum\limits_{n \geq 0}\left(  \int\limits_{\Lambda_1} f_{1,n}(\lambda_1) dA_1(\lambda_1)  \right)  X \left(  \int\limits_{\Lambda_2} f_{2,n}(\lambda_2) dA_2(\lambda_2)\right).
\end{eqnarray}

Instead the summation form provided by Eq.~\eqref{eq:integrand sum form}, the integrand function $\psi$  belongs to the \emph{integral projective tensor product} as $L^{\infty}(A_1) \acute{\otimes} L^{\infty}(A_2)$, i.e., the integrand function $\psi$ can be expressed as
\begin{eqnarray}\label{eq:integrand integral form}
\psi(\lambda_1,\lambda_2) = \int\limits_{\Xi} f_{1}(\lambda_1,x) f_{2}(\lambda_2,x) d \mu(x),
\end{eqnarray}
where $(\Xi, \mu)$ is a measure space, $f_1$ is a measurable function on $\Lambda_1 \times \Xi$, $f_2$ is a measurable function on $\Lambda_2 \times \Xi$, and 
\begin{eqnarray}
 \int\limits_{\Xi} \left\Vert   f_{1}(\cdot, x)) \right\Vert_{L^{\infty}(A_1)}
\left\Vert  f_{2}(\cdot, x)) \right\Vert_{L^{\infty}(A_2)} d \mu(x) < \infty.
\end{eqnarray}
Then, for such integrand function $\psi$, we have
\begin{eqnarray}
\int\limits_{\Lambda_1} \int\limits_{\Lambda_2} \psi(\lambda_1,\lambda_2) dA_1(\lambda_1) X dA_2(\lambda_2) \in \mathfrak{S}_1 = \int\limits_{\Xi}\left(  \int\limits_{\Lambda_1} f_{1}(\lambda_1, x) dA_1(\lambda_1)  \right)    X\left(  \int\limits_{\Lambda_2} f_{2}(\lambda_2, x) dA_2(\lambda_2)\right) d \mu(x).
\end{eqnarray}

Let $\Lambda_1, \Lambda_2$ be compact sets, and, for any $x$, suppose we have
\begin{eqnarray}\label{eq:approximation by seq}
\overline{f}_1(\lambda_1, x) &=& \lim\limits_{n \rightarrow \infty}f_{1,n}(\lambda_1, x) \nonumber \\
\overline{f}_2(\lambda_2, x) &=& \lim\limits_{n \rightarrow \infty}f_{2,n}(\lambda_2, x)
\end{eqnarray}
where $f_{1,n}$ is a measurable function on $\Lambda_1 \times \Xi$, and $f_{2,n}$ is a measurable function on $\Lambda_2 \times \Xi$. We also assume that  
\begin{eqnarray}
 \int\limits_{\Xi} \left\Vert   \overline{f}_{1}(\cdot, x)) \right\Vert_{L^{\infty}(A_1)}
\left\Vert  \overline{f}_{2}(\cdot, x)) \right\Vert_{L^{\infty}(A_2)} d \mu(x) < \infty.
\end{eqnarray}

We define the function $\psi_n(\lambda_1, \lambda_2)$ as 
\begin{eqnarray}
\psi_n(\lambda_1, \lambda_2) = \int\limits_{\Xi} f_{1,n}(\lambda_1,x) f_{2,n}(\lambda_2,x) d \mu(x).
\end{eqnarray}
Then, the limiting function of $\psi_n(\lambda_1, \lambda_2)$, denoted by $\overline{\psi}(\lambda_1, \lambda_2)$, can be expressed as 
\begin{eqnarray}
\overline{\psi}(\lambda_1, \lambda_2) = \lim\limits_{n \rightarrow \infty}\psi_n(\lambda_1, \lambda_2).
\end{eqnarray}

If functions $f_{1,n}, f_{2,n}$ are increasing with $n$, or they are dominated by some integrable functions, we have
\begin{eqnarray}\label{eq:integrand int form limit}
\overline{\psi}(\lambda_1, \lambda_2) &=& \lim\limits_{n \rightarrow \infty}\psi_n(\lambda_1, \lambda_2)
= \lim\limits_{n \rightarrow \infty}\int\limits_{\Xi} f_{1,n}(\lambda_1,x) f_{2,n}(\lambda_2,x) d \mu(x) \nonumber \\
 &=& \int\limits_{\Xi}  \lim\limits_{n \rightarrow \infty} ( f_{1,n}(\lambda_1,x) f_{2,n}(\lambda_2,x) ) d \mu(x)
\nonumber \\
 &=& \int\limits_{\Xi} \overline{f}_{1}(\lambda_1,x)\overline{f}_{2}(\lambda_2,x) d \mu(x).
\end{eqnarray}

We can summarize the above completeness arguments with Theorem on Schur multipliers given by~\cite{peller2006multiple} to have the following new theorem after taking the limit of the integrand function.
\begin{theorem}\label{thm:limit schur multiplier}
Let $\overline{\psi}$ be a measurable function of $\Lambda_1 \times \Lambda_2$, the following statements are equivalent:
\begin{enumerate}
\item $\overline{\psi} \in \Gamma(A_1, A_2)$;
\item $\overline{\psi} \in L^{\infty}(A_1) \acute{\otimes} L^{\infty}(A_2)$;
\item there exist measurable functions  $\overline{f}_{1} = \lim\limits_{n \rightarrow \infty}f_{1,n}(\lambda_1, x)$ is a measurable function on $\Lambda_1 \times \Xi$, and $\overline{f}_{2} = \lim\limits_{n \rightarrow \infty}f_{2,n}(\lambda_2, x)$ is a measurable function on $\Lambda_2 \times \Xi$ such that Eq.~\eqref{eq:integrand int form limit} holds and \\$ \int\limits_{\Xi} \left\Vert   \overline{f}_{1}(\cdot, x)) \right\Vert_{L^{\infty}(A_1)}
\left\Vert  \overline{f}_{2}(\cdot, x)) \right\Vert_{L^{\infty}(A_2)} d \mu(x) < \infty$.
\end{enumerate}
\end{theorem}

\begin{remark}
If we can find functions $
\overline{f}_{1,n}(\lambda_1) = \lim\limits_{m \rightarrow \infty}f_{1,n,m}(\lambda_1)$ and $
\overline{f}_{2,n}(\lambda_2) = \lim\limits_{m \rightarrow \infty}f_{2,n,m}(\lambda_2)$ such that 
functions $f_{1,n,m}, f_{2,n,m}$ are increasing with $m$, or they are dominated by some integrable functions,  Theorem~\ref{thm:limit schur multiplier} will still be true by replacing the statement 3 with 
\begin{eqnarray}\label{eq:integrand sum form limit}
\overline{\psi}(\lambda_1,\lambda_2) = \sum\limits_{n \geq 0} \overline{f}_{1,n}(\lambda_1) \overline{f}_{2,n}(\lambda_2),
\end{eqnarray}
and $\sum\limits_{n \geq 0}  \left\Vert \overline{f}_{1,n} \right\Vert_{L^{\infty}(A_1)}
 \left\Vert \overline{f}_{2,n} \right\Vert_{L^{\infty}(A_12} < \infty $. 
\end{remark}

If a given bivariate function can be expressed by the right-hand side of Eq,~\eqref{eq:integrand sum form limit}, this function is said constructible by the \emph{limit of projective tensor product, $\hat{\otimes}$}. On the other hand, if a given bivariate function can be expressed by the right-hand side of Eq,~\eqref{eq:integrand int form limit}, this function is said constructible by the \emph{limit of integral projective tensor product, $\acute{\otimes}$}.   


\subsection{Multiple Operator Integrals}\label{sec:MOI}

We can easily extend the definition of the limit of projective tensor product and the limit of integral projective tensor product to three or more function spaces. We have following two multiple operator integrals for the integrand function constructible by \emph{limit of projective tensor product, $\hat{\otimes}$} and by \emph{limit of integral projective tensor product, $\acute{\otimes}$}.  We use $\underline{A}_{i}^{j}$ to represent the sequence of operators $A_i, A_{i+1}, \cdots, A_j$, similarly, we apply $\underline{X}_{i}^{j}$ to represent the sequence of operators $X_i, X_{i+1}, \cdots, X_j$.

\begin{definition}\label{def:MOI sum}
Let $(\Lambda_1, A_1), (\Lambda_2, A_2), \cdots, (\Lambda_m, A_m)$ be spaces with spectral measures $\underline{A}_{1}^m$ on a HIlbert space $\mathfrak{H}$. Given a bounded measurable function $\overline{\psi}_{\hat{\otimes}}(\lambda_1, \lambda_2, \cdots, \lambda_{m})$ ($m$-variables) constructible by the \textbf{limit of projective tensor product}, and a set of $(m-1)$ operators $\underline{X}_1^{m-1}$ on Hilber space $\mathfrak{H}$, the \emph{Multiple Operator Integrals}(MOIs), denoted as $T^{\underline{A}_{1}^m}_{\overline{\psi}_{\hat{\otimes}}}(\underline{X}_{1}^{m-1})$, is defined as
\begin{eqnarray}\label{eq:MOI def}
\lefteqn{T^{\underline{A}_{1}^m}_{\overline{\psi}_{\hat{\otimes}}}(\underline{X}_1^{m-1}) \define}
\nonumber \\
&& \int\limits_{\Lambda_1}\int\limits_{\Lambda_2}\cdots \int\limits_{\Lambda_m}  \overline{\psi}_{\hat{\otimes}}(\lambda_1,\lambda_2, \cdots, \lambda_m) dA_1(\lambda_1) X_1 dA_2(\lambda_2) X_2 \cdots X_{m-1}dA_m(\lambda_m),
\end{eqnarray}
where $\overline{\psi}_{\hat{\otimes}}(\lambda_1,\lambda_2, \cdots, \lambda_m)$ can be expressed as 
\begin{eqnarray}
\overline{\psi}_{\hat{\otimes}}(\lambda_1,\lambda_2, \cdots, \lambda_m) 
&=& \sum\limits_{n \geq 0} \overline{f}_{1,n}(\lambda_1) \overline{f}_{2,n}(\lambda_2) \cdots 
\overline{f}_{m,n}(\lambda_m).
\end{eqnarray}
\end{definition}

\begin{definition}\label{def:MOI int}
Let $(\Lambda_1, A_1), (\Lambda_2, A_2), \cdots, (\Lambda_m, A_m)$ be spaces with spectral measures $\underline{A}_{1}^m$ on a HIlbert space $\mathfrak{H}$. Given a bounded measurable function $\overline{\psi}_{\acute{\otimes}}(\lambda_1, \lambda_2, \cdots, \lambda_{m})$ ($m$-variables) constructible by the \textbf{limit of integral projective tensor product}, and a set of $(m-1)$ operators $\underline{X}_{1}^{m-1}$ on Hilber space $\mathfrak{H}$, the \emph{Multiple Operator Integrals}(MOIs), denoted as $T^{\underline{A}_{1}^m}_{\overline{\psi}_{\acute{\otimes}}}(\underline{X}_1^{m-1})$, is defined as
\begin{eqnarray}\label{eq:MOI def}
T^{\underline{A}_{1}^m}_{\overline{\psi}_{\acute{\otimes}}}(\underline{X}_1^{m-1})&\define& \int\limits_{\Lambda_1}\int\limits_{\Lambda_2}\cdots \int\limits_{\Lambda_m}  \overline{\psi}_{\acute{\otimes}}(\lambda_1,\lambda_2, \cdots, \lambda_m) dA_1(\lambda_1) X_1 dA_2(\lambda_2) X_2 \cdots X_{m-1}dA_m(\lambda_m),
\end{eqnarray}
where $\overline{\psi}_{\acute{\otimes}}(\lambda_1,\lambda_2, \cdots, \lambda_m)$ can be expressed as 
\begin{eqnarray}
\overline{\psi}_{\acute{\otimes}}(\lambda_1,\lambda_2, \cdots, \lambda_m) 
&=& \int\limits_{\Xi} \overline{f}_{1}(\lambda_1,x)\overline{f}_{2}(\lambda_2,x) \cdots \overline{f}_{m}(\lambda_m,x) d \mu(x).
\end{eqnarray}
\end{definition}



Below, we will extend our previous work about DOI definition of Hermitian tensors to MOI definition based on Definitions~\ref{def:MOI sum} or~\ref{def:MOI int}. We use the symbol $\mathcal{A} \star_k \mathcal{B}$ to represent the multiplication operation between two tensors with $k$ common indices~\cite{liang2019further}. From Theorem 3.2 in~\cite{liang2019further}, every Hermitian tensor $\mathcal{H} \in  \mathbb{C}^{I_1 \times \cdots \times I_N \times I_1 \times \cdots \times I_N}$ has the following decomposition
\begin{eqnarray}\label{eq:Hermitian Eigen Decom}
\mathcal{H} &=& \sum\limits_{i=1}^{\mathbb{I}_{1}^{N}} \lambda_i \mathcal{U}_i \star_1 \mathcal{U}^{H}_i  \mbox{
~with~~$\langle \mathcal{U}_i, \mathcal{U}_i \rangle =1$ and $\langle \mathcal{U}_i, \mathcal{U}_j \rangle = 0$ for $i \neq j$,} \nonumber \\
&\define& \sum\limits_{i=1}^{\mathbb{I}_{1}^{N}} \lambda_i \mathcal{P}_{\mathcal{U}_i}
\end{eqnarray}
where $ \mathcal{U}_i \in  \mathbb{C}^{I_1 \times \cdots \times I_N \times 1}$, and the tensor $\mathcal{P}_{\mathcal{U}_i}$ is defined as $\mathcal{U}_i \star_1 \mathcal{U}^{H}_i$. The values $\lambda_i$ are named as \emph{eigevalues}. A Hermitian tensor with the decomposition shown by Eq.~\eqref{eq:Hermitian Eigen Decom} is named as \emph{eigen-decomposition}. A Hermitian tensor $\mathcal{H}$ is a positive definite (or positive semi-definite) tensor if all its eigenvalues are positive (or nonnegative).  

Let $\mathcal{A}_i \in  \mathbb{C}^{I_1 \times \cdots \times I_N \times I_1 \times \cdots \times I_N}$ for $i=1,2, \cdots, m$ be $m$ Hermitian tensors with the following eigen-decompositions:
\begin{eqnarray}
\mathcal{A}_i &=& \sum\limits_{j=1}^{\mathbb{I}_{1}^{N}} \lambda_{i,j} \mathcal{U}_{i,j} \star_1 \mathcal{U}^{H}_{i,j}
\define \sum\limits_{j=1}^{\mathbb{I}_{1}^{N}} \lambda_{i,j} \mathcal{P}_{\mathcal{U}_{i,j}}.\end{eqnarray}
We define \emph{multiple tensor integrals} (MTIs) with respect to tensors $\mathcal{A}_i$ and the function $\psi: \mathbb{R}^m \rightarrow \mathbb{R}$ which is either constructible by the \emph{limit of projective tensor product} or by the \emph{limit of integral projective tensor product}, denoted as $T^{\underline{\mathcal{A}}_1^m}_{\overline{\psi}}(\underline{\mathcal{X}}_1^{m-1})$, which can be expressed as
\begin{eqnarray}\label{eq:MTI Def}
\lefteqn{T^{\underline{\mathcal{A}}_1^m}_{\overline{\psi}}(\underline{\mathcal{X}}_1^{m-1}) =} \nonumber \\
&&
\sum\limits_{i_1=1}^{\mathbb{I}_{1}^{N}}\sum\limits_{i_2=1}^{\mathbb{I}_{1}^{N}} \cdots  \sum\limits_{i_m=1}^{\mathbb{I}_{1}^{N}} \overline{\psi}(\lambda_1, \lambda_2, \cdots, \lambda_m)  \mathcal{P}_{\mathcal{U}_{1, i_{1}}} \star_N \mathcal{X}_1 \star_N  \mathcal{P}_{\mathcal{U}_{2, i_{1}}} \star_N \mathcal{X}_2 \star_N \cdots  \star_N \mathcal{X}_{m-1} \star_N \mathcal{P}_{\mathcal{U}_{m,i_m}},
\end{eqnarray}
where $\mathcal{X}_i  \in  \mathbb{C}^{I_1 \times \cdots \times I_N \times I_1 \times \cdots \times I_N}$ for $i=1,2,\cdots, m-1$.

\section{Integrand Approximation By Linear Functions}\label{sec:Integrand Approximation By Linear Functions}

In this section, we will show that it is possible to represent any function $f(\lambda_1, \lambda_2, \cdots, \lambda_m)$ defined on a compact set of a Euclidean space $\mathbb{R}^m$ by the \emph{limit of projective tensor product} of linear functions. 

We begin with the following theorem from~\cite{chui1993realization} which says that any multivariate polynomial in $\mathbb{R}^m$ can be expressed by a polynomial of inner product forms. 
\begin{theorem}\label{thm:realization by inner product}
Let $\mathfrak{P}_k^m$ denote the set of all polynomials in $\mathbb{R}^{m}$ with polynomial degree at most $k$, and let $\bm{x}^{\bm{j}_{i,d}} = x_1^{j_{i,d,1}}x_2^{j_{i,d,2}} \cdots x_m^{j_{i,d,m}}$ such that $j_{i,d,1} + j_{i,d,2} + \cdots + j_{i,d,m} = i$  for any $d$. For any polynomial in $\mathfrak{P}_k^m$, we can express it as
\begin{eqnarray}\label{eq1:realization by inner product}
p(\bm{x}) = \sum\limits_{i=0}^k \sum\limits_{d=1}^{n_i} a_{i, d}\bm{x}^{\bm{j}_{i,d}},
\end{eqnarray}
where coefficients $a_{i, d} \in \mathbb{R}$ and $n_i = {m+i-1 \choose i}$. Note that the number $n_i$ indicates the total number of $i$-th degree polynomials in $m$ variables. Then, there exist a set of coefficients $c_{i, d} \in \mathbb{R}$ and a set of real vectors $\bm{v}_{i, d} \in \mathbb{R}^m$ such that Eq.~\eqref{eq1:realization by inner product} can be expressed as:
\begin{eqnarray}\label{eq2:realization by inner product}
p(\bm{x}) = \sum\limits_{i=0}^k \sum\limits_{d=1}^{n_i} c_{i, d}\langle \bm{x}, \bm{v}_{i, d}\rangle^i.
\end{eqnarray}
\end{theorem}

Following theorem will apply Theorem~\ref{thm:realization by inner product} to show that any multivariate function $f(\lambda_1, \lambda_2, \cdots, \lambda_m)$ defined on a compact set of a Euclidean space $\mathbb{R}^m$ can be approximated by the \emph{limit of projective tensor product} of linear functions. 
\begin{theorem}\label{thm:realization by product of linear}
Given a continous function $f(\bm{x})$ with $\bm{x} \in \mathbb{R}^m$ defined on a compact set $\mathscr{C} \in \mathbb{R}^m$ and any $\epsilon >0$, we can find the following polynomial $p(\bm{x})$ for all $\bm{x} \in \mathscr{C}$ expressed as
\begin{eqnarray}
p(\bm{x}) = \sum\limits_{i=1}^N \prod_{j=1}^i \langle \check{\bm{x}},  \bm{u}_{i, j} \rangle,
\end{eqnarray}
where $\bm{u}_{i, j} \in \mathbb{R}^{m+1}$ and $\check{\bm{x}} = [\bm{x}, 1]$; such that 
\begin{eqnarray}\label{eq1:thm:realization by product of linear}
\left\vert f(\bm{x}) - p(\bm{x})   \right\vert < \epsilon.
\end{eqnarray}
Note that $N = \sum\limits_{i=0}^{k}n_i$, where $n_i$ comes from the expression for $p(\bm{x})$ provided by Eq.~\eqref{eq1:realization by inner product}.
\end{theorem}
\textbf{Proof:}

Our first goal is to show that $p(\bm{x}) \in \mathfrak{P}_k^m$ shown by Eq.~\eqref{eq1:realization by inner product} can be expressed as
\begin{eqnarray}\label{eq2:thm:realization by product of linear}
p(\bm{x}) = \sum\limits_{i=1}^N \prod\limits_{j=1}^{i} \langle \check{\bm{x}},  \bm{u}_{i, j} \rangle ,
\end{eqnarray}
where $N = \sum\limits_{i=0}^{k} n_i$ and $\bm{u}_{i, j} \in \mathbb{R}^{m+1}$. 

From Theorem~\ref{thm:realization by inner product}, we have
\begin{eqnarray}\label{eq3:thm:realization by product of linear}
p(\bm{x}) &=& \sum\limits_{i=0}^k \sum\limits_{d=1}^{n_i} c_{i, d}\langle \bm{x}, \bm{v}_{i, d}\rangle^i \nonumber \\
&=& \sum\limits_{i=1}^{n_0}(\langle \bm{x}, \bm{0} \rangle + c_{0,i}) ( \langle \bm{x}, \bm{0} \rangle + 1 )^{i-1} \nonumber \\
&& +  \sum\limits_{i=n_0 + 1}^{n_0+ n_1}(\langle \bm{x}, c_{1,i - n_0} \bm{v}_{1,i - n_0} \rangle + 0) ( \langle \bm{x}, \bm{0} \rangle + 1 )^{i-1} \nonumber \\
&& + \cdots \nonumber \\
&& + \sum\limits_{i=n_0 + \cdots + n_{k-1} + 1}^{n_0+ \cdots + n_k} \Big[ (\langle \bm{x}, c_{k,i - (n_0+\cdots+n_{k-1})} \bm{v}_{k,i - (n_0+\cdots+n_{k-1})} \rangle + 0)  \nonumber \\
&& ~~~ \cdot ( \langle \bm{x}, \bm{v}_{k,i - (n_0+\cdots+n_{k-1})} \rangle + 0 )^{k-1}( \langle \bm{x}, \bm{0} \rangle + 1 )^{i-k} \Big].
\end{eqnarray}

On the other hand, from Eq.~\eqref{eq2:thm:realization by product of linear}, we can express $p(\bm{x})$ as 
\begin{eqnarray}\label{eq4:thm:realization by product of linear}
p(\bm{x}) &=& \sum\limits_{i=1}^N \prod\limits_{j=1}^{i} \langle \check{\bm{x}},  \bm{u}_{i, j} \rangle \nonumber \\
&=& \sum\limits_{i=1}^{n_0}\prod\limits_{j=1}^{i}\langle \check{\bm{x}}, \bm{u}_{i, j}\rangle \nonumber \\
&& +  \sum\limits_{i=n_0 + 1}^{n_0+ n_1}\prod\limits_{j=1}^{i} \langle \check{\bm{x}}, \bm{u}_{i, j} \rangle + \cdots \nonumber \\
&& + \sum\limits_{i=n_0 + \cdots + n_{k-1} + 1}^{n_0+ \cdots + n_k}\prod\limits_{j=1}^{i} \langle \check{\bm{x}}, \bm{u}_{i, j} \rangle. 
\end{eqnarray}
Because the expression in Eq.~\eqref{eq3:thm:realization by product of linear} is a special case of Eq.~\eqref{eq4:thm:realization by product of linear}, we have the expression of polynomial $p(\bm{x}) \in \mathfrak{P}_k^m$ indicated by Eq.~\eqref{eq2:realization by inner product}. Finally, this theorem is proved from generalized Weierstrass polynomial approximation Theorem~\cite{perez2006survey}.
$\hfill \Box$


\section{MOIs Properties and Random MOIs}\label{sec:MOIs Properties and Random MOIs}

In this section, we will discuss MOI properties and random MOI for the integrand function constructible by the \emph{limit of integral projective tensor product, $\acute{\otimes}$}. Same arguments will also be applied to  integrand function constructible by the \emph{limit of projective tensor product, $\hat{\otimes}$}.
 
\subsection{Properties of MOI}\label{sec:Properties of MOI}

\subsubsection{Algebraic Properties of MOI}\label{sec:Algebraic Properties of MOI}


\begin{theorem}\label{thm:MOI Algebraic Prop}
Let $(\Lambda_1, A_1), (\Lambda_2, A_2), \cdots, (\Lambda_m, A_m)$ be spaces with spectral measures $\underline{A}_{1}^m$ on a HIlbert space $\mathfrak{H}$.
\begin{roster}
\item Let $\overline{\phi}_{\acute{\otimes}}, \overline{\psi}_{\acute{\otimes}} \in L^{\infty}(A_1) \acute{\otimes} \cdots \acute{\otimes}  L^{\infty}(A_m)$ and $\alpha, \beta \in \mathbb{C}$, then we have
\begin{eqnarray}
T^{\underline{A}_1^m}_{\alpha\overline{\phi}_{\acute{\otimes}} + \beta \overline{\psi}_{\acute{\otimes}}}(\underline{X}_1^{m-1}) =  \alpha T^{\underline{A}_1^m}_{\overline{\phi}_{\acute{\otimes}}}(\underline{X}_1^{m-1}) + \beta T^{\underline{A}_1^m}_{\overline{\psi}_{\acute{\otimes}}}(\underline{X}_1^{m-1}). 
\end{eqnarray}
\item  Given $\overline{\psi}_{\acute{\otimes}, 1},  \in L^{\infty}(A_1) \acute{\otimes} \cdots \acute{\otimes}  L^{\infty}(A_k)$ and $\overline{\psi}_{\acute{\otimes}, 2},  \in L^{\infty}(A_{k+1}) \acute{\otimes} \cdots \acute{\otimes}  L^{\infty}(A_m)$, we define $\overline{\psi}_{\acute{\otimes}, 1}\oplus \overline{\psi}_{\acute{\otimes}, 2}$ as $\overline{\psi}_{\acute{\otimes}, 1}\oplus \overline{\psi}_{\acute{\otimes}, 2} \define \overline{\psi}_{\acute{\otimes}, 1}\overline{\psi}_{\acute{\otimes}, 2}$, then we have 
\begin{eqnarray}\label{eq1:thm:MOI Algebraic Prop}
\overline{\psi}_{\acute{\otimes}, 1}\oplus \overline{\psi}_{\acute{\otimes}, 2}&\in&  L^{\infty}(A_1) \acute{\otimes} \cdots \acute{\otimes}  L^{\infty}(A_m)
\end{eqnarray}
and
\begin{eqnarray}\label{eq2:thm:MOI Algebraic Prop}
T^{(\underline{A}_1^{m})}_{\overline{\psi}_{\acute{\otimes}, 1}\oplus \overline{\psi}_{\acute{\otimes}, 2}}(\underline{X}_1^{m-1}) = T^{(\underline{A}_1^{k})}_{\overline{\phi}_{\acute{\otimes},1}}(\underline{X}_1^{k-1}) X_k T^{(\underline{A}_{k+1}^{m})}_{\overline{\psi}_{\acute{\otimes},2}}(\underline{X}_{k+1}^{m-1}). 
\end{eqnarray}
\item Let $\overline{\phi}_{\acute{\otimes}} \in L^{\infty}(A_1) \acute{\otimes} \cdots \acute{\otimes}  L^{\infty}(A_m)$, and $\overline{\psi}_{\acute{\otimes},i}  \in L^{\infty}(A_{i,1}) \acute{\otimes} \cdots \acute{\otimes}  L^{\infty}(A_{i,j_i})$ for $i=1,2,\cdots, \ell$, be integrand functions composed by the $i$-th partition from $\underline{A}_{1}^m$ into $\ell$ segments such that each segment will have at least one element. For example, $(A_1, A_2, A_3, A_4, A_5)$ is partitioned into $(A_1, A_2)$, $(A_3)$, $(A_4, A_5)$. Suppose we have the following relationship between the function $\overline{\phi}_{\acute{\otimes}}$ and functions $\overline{\psi}_{\acute{\otimes},i}$ for $i=1,2,\cdots, \ell$: 
\begin{eqnarray}\label{eq3:thm:MOI Algebraic Prop}
\overline{\phi}_{\acute{\otimes}}(\lambda_1, \cdots, \lambda_m) = \prod\limits_{i=1}^{\ell} \overline{\psi}_{\acute{\otimes},i}(\underline{\lambda}_i),
\end{eqnarray}
where $\underline{\lambda}_i$ are those $j_i$ eigenvalues corresponding to spectral measures $A_{i,1}, \cdots, A_{i,j_i}$. Then, we hve
\begin{eqnarray}\label{eq4:thm:MOI Algebraic Prop}
T^{\underline{A}_1^{m}}_{\overline{\phi}_{\acute{\otimes}}}(\underline{X}_1^{m-1}) = \left(\prod\limits_{i=1}^{\ell - 1} T^{\underline{A}_{j_{i-1}+1}^{j_{i}}}_{\overline{\psi}_{\acute{\otimes},i}}(\underline{X}_{j_{i-1}+1}^{j_{i} - 1})X_{j_i}\right)T^{\underline{A}_{j_{\ell-1}+1}^{j_{\ell}}}_{\overline{\psi}_{\acute{\otimes},\ell}}(\underline{X}_{j_{\ell - 1}+1}^{j_{\ell} - 1}),
\end{eqnarray}
where $j_0 = 0$ and $j_{\ell} = m$.
\end{roster}
\end{theorem}
\textbf{Proof:}

The proof for part (a) is trivial from linearity of integration. 

For part (b), because we have
\begin{eqnarray}
\overline{\psi}_{\acute{\otimes},1}(\lambda_1,\lambda_2, \cdots, \lambda_k) 
&=& \int\limits_{\Xi_1} \overline{f}_{1}(\lambda_1,x)\overline{f}_{2}(\lambda_2,x) \cdots \overline{f}_{k}(\lambda_k,x) d \mu_1(x),
\end{eqnarray} 
and
\begin{eqnarray}
\overline{\psi}_{\acute{\otimes},2}(\lambda_{k+1},\lambda_{k+2}, \cdots, \lambda_m) 
&=& \int\limits_{\Xi_2} \overline{f}_{{k+1}}(\lambda_{k+1},x)\overline{f}_{k+2}(\lambda_{k+2},x) \cdots \overline{f}_{m}(\lambda_m,x) d \mu_2(x),
\end{eqnarray} 
then, 
\begin{eqnarray}
\overline{\psi}_{\acute{\otimes},1} \oplus \overline{\psi}_{\acute{\otimes},2} 
&=& \left( \int\limits_{\Xi_1} \overline{f}_{1}(\lambda_1,x)\overline{f}_{2}(\lambda_2,x) \cdots \overline{f}_{k}(\lambda_k,x) d \mu_1(x)\right) \nonumber \\
&& \cdot 
\left( \int\limits_{\Xi_2} \overline{f}_{{k+1}}(\lambda_{k+1},x)\overline{f}_{k+2}(\lambda_{k+2},x) \cdots \overline{f}_{m}(\lambda_m,x) d \mu_2(x)\right)
 \nonumber \\
&=&  \int\limits_{\Xi_1 \times \Xi_2} \overline{f}_{1}(\lambda_1,x) \cdots \overline{f}_{k}(\lambda_k,x) \overline{f}_{{k+1}}(\lambda_{k+1},x) \cdots \overline{f}_{m}(\lambda_m,x)   d \mu_1(x) d \mu_2(x).
\end{eqnarray} 
Eq.~\eqref{eq1:thm:MOI Algebraic Prop} is true from the product measure space. 

Since we have
\begin{eqnarray}
T^{\underline{A}_{1}^{m}}_{\overline{\psi}_{\acute{\otimes}, 1}\oplus \overline{\psi}_{\acute{\otimes}, 2}}(\underline{X}_{1}^{m-1}) &=& \int\limits_{\Lambda_1}\int\limits_{\Lambda_2}\cdots \int\limits_{\Lambda_m}  \overline{\psi}_{\acute{\otimes}, 1}\overline{\psi}_{\acute{\otimes}, 2} dA_1(\lambda_1) X_1 dA_2(\lambda_2) X_2 \cdots X_{m-1} dA_m(\lambda_m) \nonumber \\
&=& \left(\int\limits_{\Lambda_1}\cdots \int\limits_{\Lambda_k}  \overline{\psi}_{\acute{\otimes}, 1} dA_1(\lambda_1) X_1 d\cdots X_{k-1} dA_k(\lambda_k)\right) X_k \nonumber \\
& & \int\limits_{\Lambda_{k+1}}\cdots \int\limits_{\Lambda_m}  \overline{\psi}_{\acute{\otimes}, 2} dA_{k+1}(\lambda_{k+1}) X_{k+1} \cdots X_{m-1} dA_m(\lambda_m) \nonumber \\
&=&  T^{(\underline{A}_1^{k})}_{\overline{\phi}_{\acute{\otimes}}, 1}(\underline{X}_1^{k-1}) X_k T^{(\underline{A}_{k+1}^{m})}_{\overline{\psi}_{\acute{\otimes}},2}(\underline{X}_{k+1}^{m-1}). 
\end{eqnarray}
Therefore, we have Eq.~\eqref{eq2:thm:MOI Algebraic Prop}.

For part (c), we have
\begin{eqnarray}
T^{\underline{A}_1^{m}}_{\overline{\phi}_{\acute{\otimes}}}(\underline{A}_1^{m-1})&=&   \int\limits_{\Lambda_1}\int\limits_{\Lambda_2}\cdots \int\limits_{\Lambda_m}  \overline{\psi}_{\acute{\otimes}} dA_1(\lambda_1) X_1 dA_2(\lambda_2) X_2 \cdots X_{m-1} dA_m(\lambda_m) \nonumber \\
&=_1&   \int\limits_{\Lambda_1}\int\limits_{\Lambda_2}\cdots \int\limits_{\Lambda_m}\left( \prod\limits_{i=1}^{\ell} \overline{\psi}_{\acute{\otimes},i}(\underline{\lambda}_i)  \right) dA_1(\lambda_1) X_1 dA_2(\lambda_2) X_2 \cdots X_{m-1} dA_m(\lambda_m) \nonumber \\
&=&   \underbrace{\left(\int\limits_{\Lambda_1}\cdots \int\limits_{\Lambda_{j_1}} \overline{\psi}_{\acute{\otimes},1}(\underline{\lambda}_1)dA_1(\lambda_1) X_1 \cdots X_{j_{1}-1} dA_{j_1}(\lambda_{j_1})\right)}_{T^{\underline{A}_1^{j_1}}_{\overline{\psi}_{\acute{\otimes},1}}(\underline{X}_1^{j_{1} - 1})} X_{j_1} \nonumber \\
& & \cdot \underbrace{\left(\int\limits_{\Lambda_{j_1+1}}\cdots \int\limits_{\Lambda_{j_2}} \overline{\psi}_{\acute{\otimes},2}(\underline{\lambda}_2)dA_{j_1 + 1}(\lambda_{j_1 +1}) X_{j_1 + 1} \cdots X_{j_2 -1} dA_{j_2}(\lambda_{j_2})\right)}_{T^{\underline{A}_{j_{1} + 1}^{j_2}}_{\overline{\psi}_{\acute{\otimes},2}}(\underline{X}_{j_{1}+1}^{j_{2} - 1})}X_{j_2} \nonumber \\
&& \cdots  \underbrace{\left(\int\limits_{\Lambda_{j_{\ell-1}+1}}\cdots \int\limits_{\Lambda_{j_{\ell}}} \overline{\psi}_{\acute{\otimes},\ell}(\underline{\lambda}_\ell)dA_{j_{\ell - 1} + 1}(\lambda_{j_{\ell -1} +1}) X_{j_{\ell-1} + 1} \cdots X_{j_{\ell}-1} dA_{j_{\ell}}(\lambda_{j_{\ell}})\right)}_{T^{\underline{A}_{j_{\ell-1}}^{j_{\ell}}}_{\overline{\psi}_{\acute{\otimes},\ell}}(\underline{X}_{j_{\ell - 1}+1}^{j_{\ell} - 1})} \nonumber \\
&=&  \left(\prod\limits_{i=1}^{\ell - 1} T^{\underline{A}_{j_{i-1}+1}^{j_{i}}}_{\overline{\psi}_{\acute{\otimes},i}}(\underline{X}_{j_{i-1}+1}^{j_{i} - 1})X_{j_i}\right)T^{\underline{A}_{j_{\ell-1}+1}^{j_{\ell}}}_{\overline{\psi}_{\acute{\otimes},\ell}}(\underline{X}_{j_{\ell - 1}+1}^{j_{\ell} - 1}),
\end{eqnarray}
where we apply Eq.~\eqref{eq3:thm:MOI Algebraic Prop} to $=_1$ and set $j_0 = 1$. Note that $j_\ell = m$.
$\hfill \Box$


\subsubsection{Norm Estimations for MOIs}\label{sec:Norm Estimations for MOI}


Following lemma will provide two norm estimations of MOIs. 
\begin{lemma}\label{lma:norm estimation}
According to Definition~\ref{def:MOI int}, we have the following norm estimations for Eq.~\eqref{eq:MOI def}:
\begin{roster}
\item  Given a bounded measurable function $\overline{\psi}_{\acute{\otimes}}(\lambda_1, \lambda_2, \cdots, \lambda_{m})$ ($m$-variables) constructible by \textbf{limit of integral projective tensor product}, and a set of $(m-1)$ bounded linear operators $X_1, X_2, \cdots, X_{m-1}$ on Hilber space $\mathfrak{H}$, we have
\begin{eqnarray}\label{eq1:lma:norm estimation}
\lefteqn{\left\Vert \int\limits_{\Lambda_1}\int\limits_{\Lambda_2}\cdots \int\limits_{\Lambda_m}  \overline{\psi}_{\acute{\otimes}}(\lambda_1,\lambda_2, \cdots, \lambda_m) dA_1(\lambda_1) X_1 dA_2(\lambda_2) X_2 \cdots X_{m-1}dA_m(\lambda_m) \right\Vert } \nonumber \\
&\leq&  \left\Vert \overline{\psi}_{\acute{\otimes}}(\lambda_1, \lambda_2, \cdots, \lambda_{m}) \right\Vert \prod\limits_{i=1}^{m-1} \left\Vert X_i \right\Vert. ~~~~~~~~~~~~~~~~~~~~~~~~~~~~~~~~~~~~~~~~~~~~~~~~~~~~~~~~~~~~~~~~~~~~~~~~
\end{eqnarray}
\item Given a bounded measurable function $\overline{\psi}_{\acute{\otimes}}(\lambda_1, \lambda_2, \cdots, \lambda_{m})$ ($m$-variables) constructible by \textbf{limit of integral projective tensor product}, and a set of $(m-1)$ bounded linear operators $X_1, X_2, \cdots, X_{m-1}$ on Hilber space $\mathfrak{H}$ such that $X_i \in \mathscr{S}_{p_i}$, where $\mathscr{S}_{p_i}$ is a Schatten–von Neumann class with parameter $p_i$, we also require that $\sum\limits_{i=1}^m \frac{1}{p_i} \leq 1$. Then, we have 
\begin{eqnarray}\label{eq2:lma:norm estimation}
\lefteqn{\left\Vert \int\limits_{\Lambda_1}\int\limits_{\Lambda_2}\cdots \int\limits_{\Lambda_m}  \overline{\psi}_{\acute{\otimes}}(\lambda_1,\lambda_2, \cdots, \lambda_m) dA_1(\lambda_1) X_1 dA_2(\lambda_2) X_2 \cdots X_{m-1}dA_m(\lambda_m) \right\Vert_{\mathscr{S}_{q}}}\nonumber \\
&\leq & \left\Vert \overline{\psi}_{\acute{\otimes}}(\lambda_1, \lambda_2, \cdots, \lambda_{m}) \right\Vert \prod\limits_{i=1}^{m-1} \left\Vert X_i \right\Vert_{\mathscr{S}_{p_i}},~~~~~~~~~~~~~~~~~~~~~~~~~~~~~~~~~~~~~~~~~~~~~~~~~~~~~~~~~~~~~~~~~~~~~
\end{eqnarray}
where $\frac{1}{q} = 1 -  \sum\limits_{i=1}^m \frac{1}{p_i}$.
\end{roster}
\end{lemma}
\textbf{Proof:}
By applying triangle inequality for integrals and H\"{o}lder’s inequality to Eq.~\eqref{eq:MOI def}.
$\hfill \Box$

\subsubsection{Perturbation Formula for MOI}\label{sec:Perturbation Formula for MOI}


Before presenting our perturbation formula for MOI, we have to define divided differences for function $\overline{\psi}_{\acute{\otimes}}(\lambda)$, which is a function constructible by \emph{limit of integral projective tensor product, $\acute{\otimes}$}. We assume this function can have $n$ times continuously
differentiable on $\mathbb{R}$. The divided difference $\overline{\psi}^{[n]}_{\acute{\otimes}}$ of order $n$ with respect to the $j$-th variable is defined by 
\begin{eqnarray}\label{eq:div diff def}
\lefteqn{\overline{\psi}^{[n+1]}_{\acute{\otimes}}(\lambda_0,\cdots, \lambda_{n+1})  \define} \nonumber \\
&& \frac{\overline{\psi}^{[n]}_{\acute{\otimes}}(\lambda_0, \cdots,\lambda_{j-1}, \lambda_{j+1}, \cdots, \lambda_{n+1}) -  \overline{\psi}^{[n]}_{\acute{\otimes}}(\lambda_0, \cdots,\lambda_{j-2}, \lambda_j, \cdots, \lambda_{n+1})}{\lambda_{j}- \lambda_{j-1}}.
\end{eqnarray}

We have the following perturbation formula for MOIs. 
\begin{lemma}\label{lma:MOI perturbation}
Let $\overline{\psi}_{\acute{\otimes}}$ be $n$ times differentiable on $\mathbb{R}$ and $C,D, A_1, \cdots, A_{n-1}$ be operators. Then, we have
\begin{eqnarray}
T^{\underline{A}_{1}^{j-1}, C,\underline{A}_{j}^{m}}_{\overline{\psi}^{[m]}_{\acute{\otimes}}}(\underline{X}_{1}^{m}) - T^{\underline{A}_{1}^{j-1}, D,\underline{A}_{j}^{m}}_{\overline{\psi}^{[m]}_{\acute{\otimes}}}(\underline{X}_{1}^{m}) = T^{\underline{A}_{1}^{j-1}, C, D,\underline{A}_{j}^{m}}_{\overline{\psi}^{[m+1]}_{\acute{\otimes}}}(\underline{X}_{1}^{j-1}, C-D,\underline{X}_{j}^{m}),
\end{eqnarray}
where $m \leq n-1$.
\end{lemma}
\textbf{Proof:}
For $1 \leq j \leq m+1$, we define 
\begin{eqnarray}
f_j(\lambda_0, \cdots, \lambda_{m+1})&\define& \lambda_j \overline{\psi}^{[m+1]}_{\acute{\otimes}}( \lambda_0, \cdots, \lambda_{m+1}), \nonumber \\
g_j(\lambda_0, \cdots, \lambda_{m+1})&\define& \overline{\psi}^{[m]}_{\acute{\otimes}}( \lambda_0, \cdots, \lambda_{j-1}, \lambda_{j+1}, \cdots, \lambda_{m+1}).
\end{eqnarray}
From Eq.~\eqref{eq:div diff def}, we have 
\begin{eqnarray}\label{eq1:lma:MOI perturbation}
f_j - f_{j-1} = g_j - g_{j-1}.
\end{eqnarray}

Then, we have
\begin{eqnarray}
\lefteqn{T^{\underline{A}_{1}^{j-1}, C, D,\underline{A}_{j}^{m}}_{\overline{\psi}^{[m+1]}_{\acute{\otimes}}}(\underline{X}_{1}^{j-1}, C-D,\underline{X}_{j}^{m})} \nonumber \\
&=& 
T^{\underline{A}_{1}^{j-1}, C, D,\underline{A}_{j}^{m}}_{\overline{\psi}^{[m+1]}_{\acute{\otimes}}}(\underline{X}_{1}^{j-1}, C,\underline{X}_{j}^{m})- 
T^{\underline{A}_{1}^{j-1}, C, D,\underline{A}_{j}^{m}}_{\overline{\psi}^{[m+1]}_{\acute{\otimes}}}(\underline{X}_{1}^{j-1}, D,\underline{X}_{j}^{m}) \nonumber \\
&=& T^{\underline{A}_{1}^{j-1}, C, D,\underline{A}_{j}^{m}}_{f_{j-1} }(\underline{X}_{1}^{j-1}, Id,\underline{X}_{j}^{m}) - 
T^{\underline{A}_{1}^{j-1}, C, D,\underline{A}_{j}^{m}}_{f_j }(\underline{X}_{1}^{j-1}, Id,\underline{X}_{j}^{m}) \nonumber \\
&=&T^{\underline{A}_{1}^{j-1}, C, D,\underline{A}_{j}^{m}}_{f_{j-1} -f_j}(\underline{X}_{1}^{j-1},Id ,\underline{X}_{j}^{m})  \nonumber \\
&=_1&T^{\underline{A}_{1}^{j-1}, C, D,\underline{A}_{j}^{m}}_{g_{j-1} -g_j}(\underline{X}_{1}^{j-1},Id ,\underline{X}_{j}^{m}) \nonumber \\
&=& T^{\underline{A}_{1}^{j-1}, C, D,\underline{A}_{j}^{m}}_{g_{j-1} }(\underline{X}_{1}^{j-1}, Id,\underline{X}_{j}^{m}) - 
T^{\underline{A}_{1}^{j-1}, C, D,\underline{A}_{j}^{m}}_{g_j}(\underline{X}_{1}^{j-1}, Id,\underline{X}_{j}^{m}) \nonumber \\
&=& T^{\underline{A}_{1}^{j-1}, C, \underline{A}_{j}^{m}}_{\overline{\psi}^{[m]}_{\acute{\otimes}}}(\underline{X}_{1}^{m}) - 
T^{\underline{A}_{1}^{j-1}, D, \underline{A}_{j}^{m}}_{\overline{\psi}^{[m]}_{\acute{\otimes}}}(\underline{X}_{1}^{m}),
\end{eqnarray}
where $Id$ is the identity operator and $=_1$ comes from Eq.~\eqref{eq1:lma:MOI perturbation}.
$\hfill \Box$

\subsubsection{Continuity for MOIs}\label{sec:Continuity for MOI}


In the next lemma, we establish continuity of a MOI defined by Definition~\ref{def:MOI int}.
\begin{lemma}\label{lma:MOI continuity}
Let $\overline{\psi}_{\acute{\otimes}}$ be $n$ times differentiable on $\mathbb{R}$ and $A_i^{(m)} \rightarrow A_i$ as $m \rightarrow \infty$, i.e.,  $\left\Vert A_i^{(m)} - A_i \right\Vert \rightarrow 0$ as $m \rightarrow \infty$. We use $\underline{A^{(m)}}_{i}^{j}$ to represent the sequence $A^{(m)}_i, A^{(m)}_{i+1}, \cdots, A^{(m)}_j$. Then, we have
\begin{eqnarray}
T^{\underline{A^{(m)}}_{1}^{n+1}}_{\overline{\psi}^{[n]}_{\acute{\otimes}}}(\underline{X}_{1}^{n}) \rightarrow T^{\underline{A}_{1}^{n+1}}_{\overline{\psi}^{[n]}_{\acute{\otimes}}}(\underline{X}_{1}^{n}),
\end{eqnarray}
where $X_1, X_2, \cdots, X_n$ are bounded operators by assuming that the norm for the function $\overline{\psi}^{[n+1]} _{\acute{\otimes}}$ and norms for the oprators $X_1, \cdots, X_n$ are bounded.
\end{lemma}
\textbf{Proof:}
Because we have
\begin{eqnarray}\label{eq1:lma:MOI continuity}
\left\Vert  T^{\underline{A^{(m)}}_{1}^{n+1}}_{\overline{\psi}^{[n]}_{\acute{\otimes}}}(\underline{X}_{1}^{n}) - T^{\underline{A}_{1}^{n+1}}_{\overline{\psi}^{[n]}_{\acute{\otimes}}}(\underline{X}_{1}^{n}) \right\Vert &=_1& \left\Vert  \sum\limits_{i=1}^{n+1} \left(T^{\underline{A}_{1}^{i-1}, \underline{A^{(m)}}_{i}^{n+1}}_{\overline{\psi}^{[n]}_{\acute{\otimes}}}(\underline{X}_{1}^{n}) - T^{\underline{A}_{1}^{i}, \underline{A^{(m)}}_{i+1}^{n+1}}_{\overline{\psi}^{[n]}_{\acute{\otimes}}}(\underline{X}_{1}^{n}) \right) \right\Vert \nonumber \\
&=_2& \left\Vert  \sum\limits_{i=1}^{n+1} T^{\underline{A}_{1}^{i-1}, A^{(m)}_i,  A_i, \underline{A^{(m)}}_{i+1}^{n+1}}_{\overline{\psi}^{[n+1]}_{\acute{\otimes}}}(\underline{X}_{1}^{i-1}, A^{(m)}_i - A_i ,\underline{X}_{i}^{n})\right\Vert
\end{eqnarray}
where $=_1$ comes from telescoping summation and $=_2$ comes from Lemma~\ref{lma:MOI perturbation}

Applying norm triangle inequality and Lemma~\ref{lma:norm estimation}, we can upper bound Eq.~\eqref{eq1:lma:MOI continuity} as 
\begin{eqnarray}\label{eq2:lma:MOI continuity}
\left\Vert  T^{\underline{A^{(m)}}_{1}^{n+1}}_{\overline{\psi}^{[n]}_{\acute{\otimes}}}(\underline{X}_{1}^{n})  - T^{\underline{A}_{1}^{n+1}}_{\overline{\psi}^{[n]}_{\acute{\otimes}}}(\underline{X}_{1}^{n})  \right\Vert &\leq& \left\Vert \overline{\psi}^{[n+1]}_{\acute{\otimes}} \right\Vert
\sum\limits_{i=1}^{n+1}\left( \left\Vert  A^{(m)}_i - A_i \right\Vert \prod\limits_{j=1}^{n} \left\Vert X_j \right\Vert\right).  
\end{eqnarray}
This theorem is provied when $m \rightarrow \infty$.
$\hfill \Box$

\subsection{Random MOIs}\label{sec:Random MOI}

From the MOI definition provided by Definition~\ref{def:MOI sum} and Definition~\ref{def:MOI int}, the random MOI considered in this work is to assume random operators $A_i$ having spectrum decomposition as $\int\limits_{\Lambda_{i}}\lambda_i dA_i(\lambda_i)$, where $i=1,2, \cdots, m$. The randomness of the operator $A_i$ is  determined by random variables $\lambda_i$ and random unitary operators $dA_i(\lambda_i)$ via Haar measure. The remaining parameters like the integrand function$\overline{\psi}_{\hat{\otimes}} (or \overline{\psi}_{\acute{\otimes}}$, and operators $X_1, \cdots, X_m-1$ are assumed to be deterministic. The purpose of this section is to derive tail bounds for random MOI norms and establish the continuity property of random MOI. 


\begin{theorem}\label{thm:tail bound based on norm estimation}
We have the following tail bound for the norm estimations:
\begin{roster}
\item  Given a random MOI, $T^{\underline{A}_{1}^{m+1}}_{\overline{\psi}^{[n]}_{\acute{\otimes}}}(\underline{X}_{1}^{m})$,  we have
\begin{eqnarray}\label{eq1:thm:tail bound based on norm estimation}
\mathrm{Pr}\left( \left\Vert T^{\underline{A}_{1}^{m+1}}_{\overline{\psi}^{[m]}_{\acute{\otimes}}}(\underline{X}_{1}^{m}) \right\Vert \geq \theta \right) 
&\leq&   \frac{\prod\limits_{i=1}^{m-1} \left\Vert X_i \right\Vert}{\theta}\mathbb{E}\left[\left\Vert \overline{\psi}_{\acute{\otimes}}(\lambda_1, \lambda_2, \cdots, \lambda_{m})\right\Vert\right], 
\end{eqnarray}
where $\theta > 0$. 
\item  Given a random MOI, $T^{\underline{A}_{1}^{m+1}}_{\overline{\psi}^{[m]}_{\acute{\otimes}}}(\underline{X}_{1}^{m})$,  we have
\begin{eqnarray}\label{eq2:thm:tail bound based on norm estimation}
\mathrm{Pr}\left( \left\Vert T^{\underline{A}_{1}^{m+1}}_{\overline{\psi}^{[m]}_{\acute{\otimes}}}(\underline{X}_{1}^{m}) \right\Vert_{\mathscr{S}_{q}} \geq \theta \right) 
&\leq&  \frac{\prod\limits_{i=1}^{m-1} \left\Vert X_i \right\Vert_{\mathscr{S}_{p_i}}}{\theta} \mathbb{E}\left[\left\Vert \overline{\psi}_{\acute{\otimes}}(\lambda_1, \lambda_2, \cdots, \lambda_{m}), \right\Vert\right],
\end{eqnarray}
where $\frac{1}{q} = 1 -  \sum\limits_{i=1}^m \frac{1}{p_i}$ and $\theta > 0$. 
\end{roster}
\end{theorem}
\textbf{Proof:}

For part (a), since we have 
\begin{eqnarray}
\mathrm{Pr}\left(\left\Vert T^{\underline{A}_{1}^{m+1}}_{\overline{\psi}^{[m]}_{\acute{\otimes}}}(\underline{X}_{1}^{m}) \right\Vert \geq \theta \right) &\leq_1&  \mathrm{Pr}\left( \left\Vert \overline{\psi}_{\acute{\otimes}}(\lambda_1, \lambda_2, \cdots, \lambda_{m}) \right\Vert \prod\limits_{i=1}^{m-1} \left\Vert X_i \right\Vert \geq \theta \right) \nonumber \\
&=&   \mathrm{Pr}\left( \left\Vert \overline{\psi}_{\acute{\otimes}}(\lambda_1, \lambda_2, \cdots, \lambda_{m}) \right\Vert \geq \frac{\theta}{\prod\limits_{i=1}^{m-1} \left\Vert X_i \right\Vert} \right) \nonumber \\ 
&\leq_2&  \frac{\prod\limits_{i=1}^{m-1} \left\Vert X_i \right\Vert}{\theta}\mathbb{E}\left[\left\Vert \overline{\psi}_{\acute{\otimes}}(\lambda_1, \lambda_2, \cdots, \lambda_{m})\right\Vert\right], 
\end{eqnarray}
where $\leq_1$ is due to part (a) of Lemma~\ref{lma:norm estimation}, and $\leq_2$ comes from Markov inequality.

For part (b), since we have 
\begin{eqnarray}
\mathrm{Pr}\left(\left\Vert T^{\underline{A}_{1}^{m+1}}_{\overline{\psi}^{[m]}_{\acute{\otimes}}}(\underline{X}_{1}^{m}) \right\Vert_{\mathscr{S}_{q}} \geq \theta \right) &\leq_1&  \mathrm{Pr}\left( \left\Vert \overline{\psi}_{\acute{\otimes}}(\lambda_1, \lambda_2, \cdots, \lambda_{m}) \right\Vert \prod\limits_{i=1}^{m-1} \left\Vert X_i \right\Vert_{\mathscr{S}_{p_i}} \geq \theta \right) \nonumber \\
&=&   \mathrm{Pr}\left( \left\Vert \overline{\psi}_{\acute{\otimes}}(\lambda_1, \lambda_2, \cdots, \lambda_{m}) \right\Vert \geq \frac{\theta}{\prod\limits_{i=1}^{m-1} \left\Vert X_i \right\Vert_{\mathscr{S}_{p_i}}} \right) \nonumber \\ 
&\leq_2&  \frac{\prod\limits_{i=1}^{m-1} \left\Vert X_i \right\Vert_{\mathscr{S}_{p_i}}}{\theta}\mathbb{E}\left[\left\Vert \overline{\psi}_{\acute{\otimes}}(\lambda_1, \lambda_2, \cdots, \lambda_{m})\right\Vert\right], 
\end{eqnarray}
where $\leq_1$ is due to part (b) of Lemma~\ref{lma:norm estimation}, and $\leq_2$ comes from Markov inequality.
$\hfill \Box$



In the remainder of this section, we will establish continuity of random MOI. We need the following definition about the convergence in mean for random operators. 

\begin{definition}\label{def:conv in mean}
We say that a sequence of random operator $X^{(m)}$ converges in the $r$-th mean towards the random operator $X$ with respect to the operator norm $\left\Vert \cdot \right\Vert$, if we have
\begin{eqnarray}
\mathbb{E}\left( \left\Vert X^{(m)} \right\Vert^{r} \right)~~~\mbox{exists,}
\end{eqnarray}
and
\begin{eqnarray}
\mathbb{E}\left( \left\Vert X \right\Vert^{r} \right)~~~ \mbox{exists,}
\end{eqnarray}
and
\begin{eqnarray}
\lim\limits_{m \rightarrow \infty}\mathbb{E}\left( \left\Vert X^{(m)}  - X \right\Vert^{r} \right) = 0. 
\end{eqnarray}
We adopt the notatation $X^{(m)} \xrightarrow[]{r} X$ to represent that random operators $X^{(m)}$ converges in the $r$-th mean to the random operator $X$ with respect to the norm $\left\Vert \cdot \right\Vert$.
\end{definition}

Following theorem is about the continuity in the $r$-th mean of random MOI . 

\begin{theorem}\label{thm:continuity of random MOI}
For $i=1, 2, \cdots, n+1$, we have random self-adjoint operators $A_i^{(m)}, A_i$ such that
\begin{eqnarray}\label{eq1:thm:continuity of random MOI}
A_i^{(m)} \xrightarrow[]{r} A_i ~~\mbox{as $m \rightarrow \infty$},
\end{eqnarray}
where $1 \leq r < \infty$. Moreover, the norm for the real valued function $\overline{\psi}^{[n+1]}_{\acute{\otimes}}$ and the norm for $X_j$ for $j=1,2,\cdots, n$, are assumed bounded. Then, we have  
\begin{eqnarray}
T^{\underline{A^{(m)}}_{1}^{n+1}}_{\overline{\psi}^{[n]}_{\acute{\otimes}}}(\underline{X}_{1}^{n}) \xrightarrow[]{r}
T^{\underline{A}_{1}^{n+1}}_{\overline{\psi}^{[n]}_{\acute{\otimes}}}(\underline{X}_{1}^{n}).
\end{eqnarray}
\end{theorem}
\textbf{Proof:}

From the proof in Lemma~\ref{lma:MOI continuity}, we have 
\begin{eqnarray}\label{eq2:thm:continuity of random MOI}
\left\Vert  T^{\underline{A^{(m)}}_{1}^{n+1}}_{\overline{\psi}^{[n]}_{\acute{\otimes}}}(\underline{X}_{1}^{n}) - T^{\underline{A}_{1}^{n+1}}_{\overline{\psi}^{[n]}_{\acute{\otimes}}}(\underline{X}_{1}^{n}) \right\Vert &\leq& \left\Vert \overline{\psi}^{[n+1]}_{\acute{\otimes}} \right\Vert
\sum\limits_{i=1}^{n+1}\left( \left\Vert  A^{(m)}_i - A_i \right\Vert \prod\limits_{j=1}^{n} \left\Vert X_j \right\Vert\right).  
\end{eqnarray}

By raising the power $r$ and taking the expectation at the both sides of the inequality provided by Eq.~\eqref{eq2:thm:continuity of random MOI}, we have proved this theorem by conditions given by Eq.~\eqref{eq1:thm:continuity of random MOI} and the following inequality:
\begin{eqnarray}
\left(\sum\limits_{i=1}^{n+1}a_i\right)^r \leq (n+1)^{r-1}\left(\sum\limits_{i=1}^{n+1}a_i^r\right)~~\mbox{given $a_i \geq 0$}.
\end{eqnarray}
$\hfill \Box$

\section{Random MOIs Applications}\label{sec:Random MOI Applications}

In this section, we will apply random MOIs to obtain several tail bounds related to several applications of MOIs. We will focus on functions constructible by the \emph{limit of integral projective tensor product, $\acute{\otimes}$}. Same arguments will also be applied to functions constructible by the \emph{limit of projective tensor product, $\hat{\otimes}$}.

\subsection{Tail Bound for Higher Random Operator Derivative}\label{sec:Tail Bound for Higher Random Operator Derivative}


In this section, we will derive tail bounds for higher random operato derivatives. We begin with a lemma to express the first derivative of an operator-valued function by double operator integrals.  

\begin{lemma}\label{lma:1-th derivative by DOIs}
Let $X(t)$ be self-adjoint operators indexed by $t$ with bounded norm, i.e., $\left\Vert X(t) \right\Vert < \infty$, and$\overline{\psi}_{\acute{\otimes}}$ be differentiable on $\mathbb{R}$. Then, we have
\begin{eqnarray}\label{eq1:lma:1-th derivative by DOIs}
\frac{d\overline{\psi}_{\acute{\otimes}}(X(t))}{dt}&=&T^{X(t),X(t)}_{\overline{\psi}_{\acute{\otimes}}^{[1]}}\left(\frac{d X(t)}{dt}\right).
\end{eqnarray}
\end{lemma}
\textbf{Proof:}
Suppose the function $\overline{\psi}_{\acute{\otimes}}(x)$ has the polynomial form as $\overline{\psi}_{\acute{\otimes}}(x) = x^m$, by the chain rule of dervative, we have 
\begin{eqnarray}\label{eq2:lma:1-th derivative by DOIs}
\frac{d\overline{\psi}_{\acute{\otimes}}(X(t))}{dt}&=& \sum\limits_{i=0}^{m-1}X^{i}(t) \frac{d X(t)}{dt} X^{m-i-1}(t) .
\end{eqnarray}
Since the operator $X(t)$ can be represented by its spectrum as
\begin{eqnarray}\label{eq3:lma:1-th derivative by DOIs}
X(t) = \int\limits_{\Lambda_t} \lambda dA_t(\lambda),
\end{eqnarray}
where $(\Lambda_t, A_t)$ is the spectral measure at the instant $t$;by the spectral theorem, we can express Eq.~\eqref{eq2:lma:1-th derivative by DOIs} as 
\begin{eqnarray}\label{eq4:lma:1-th derivative by DOIs}
\frac{d\overline{\psi}_{\acute{\otimes}}(X(t))}{dt}&=& \sum\limits_{i=0}^{m-1} \int\limits_{\Lambda_t} \lambda_1^i dA_t(\lambda_1) \frac{d X(t)}{dt} \int\limits_{\Lambda_t} \lambda_2^{m-i-1} dA_t(\lambda_2) \nonumber \\
&=&  \int\limits_{\Lambda_t}  \int\limits_{\Lambda_t}  \left(\sum\limits_{i=0}^{m-1} \lambda_1^i  \lambda_2^{m-i-1}\right) dA_t(\lambda_1) \frac{d X(t)}{dt} dA_t(\lambda_2) \nonumber \\
&=&  \int\limits_{\Lambda_t}  \int\limits_{\Lambda_t} \underbracket{\frac{\lambda_1^m - \lambda_2^m}{\lambda_1 - \lambda_2}}_{(x^m)^{[1]}}dA_t(\lambda_1) \frac{d X(t)}{dt} dA_t(\lambda_2).
\end{eqnarray}
From MOIs defined by Eq.~\eqref{def:MOI int} as $m=2$ (DOIs), we have Eq.~\eqref{eq1:lma:1-th derivative by DOIs} for $\overline{\psi}_{\acute{\otimes}}(x) = x^m$. By linearity, this lemma is also true when $\overline{\psi}_{\acute{\otimes}}(x)$ is a polynomial function. 

Suppose we can approximate $\overline{\psi}_{\acute{\otimes}}$ by a sequence of polynomials $p_n(x)$, i.e., $\left\Vert \overline{\psi}_{\acute{\otimes}}-p_n \right\Vert \rightarrow 0$ as $n \rightarrow \infty$. Then, we have
\begin{eqnarray}
\left\Vert T^{X(t),X(t)}_{\overline{\psi}_{\acute{\otimes}}^{[1]}}\left(\frac{d X(t)}{dt}\right) - 
T^{X(t),X(t)}_{p_n^{[1]}}\left(\frac{d X(t)}{dt}\right)  \right\Vert \leq Const \left\Vert \overline{\psi}_{\acute{\otimes}}-p_n \right\Vert \left\Vert X(t) \right\Vert.
\end{eqnarray}
This lemma is proved due to we can approximate $\overline{\psi}_{\acute{\otimes}}$ by polynomials. 
$\hfill \Box$

From Lemma~\ref{lma:1-th derivative by DOIs}, we are ready to have the following theorem about the tail bound for $\left\Vert \frac{d\overline{\psi}_{\acute{\otimes}}(X(t))}{dt} \right\Vert$.
\begin{theorem}\label{thm:tail bound for the first der of t}
Let $X(t)$ be self-adjoint operators indexed by $t$ with bounded norm, i.e., $\left\Vert X(t) \right\Vert < \infty$, and$\overline{\psi}_{\acute{\otimes}}$ be differentiable on $\mathbb{R}$. We also assume that $\left\Vert \frac{d X(t)}{dt} \right\Vert < \Upsilon_X$. Then, 
\begin{eqnarray}\label{eq1:thm:tail bound for the first der of t}
\mathrm{Pr}\left( \left\Vert \frac{d\overline{\psi}_{\acute{\otimes}}(X(t))}{dt} \right\Vert > \theta \right) & \leq & \frac{\Upsilon_X}{\theta} \mathbb{E}\left[ \left\Vert \overline{\psi}^{[1]}_{\acute{\otimes}} \right\Vert  \right].
\end{eqnarray}
\end{theorem}
\textbf{Proof:}
Because we have
\begin{eqnarray}
\mathrm{Pr}\left( \left\Vert \frac{d\overline{\psi}_{\acute{\otimes}}(X(t))}{dt} \right\Vert > \theta \right) & = & \mathrm{Pr}\left( \left\Vert T^{X(t),X(t)}_{\overline{\psi}_{\acute{\otimes}}^{[1]}}\left(\frac{d X(t)}{dt}\right) \right\Vert > \theta \right) \nonumber \\
& \leq_1 & \mathrm{Pr}\left( \left\Vert \overline{\psi}^{[1]}_{\acute{\otimes}} \right\Vert \left\Vert \frac{d X(t)}{dt} \right\Vert > \theta \right) \nonumber \\
& = & \mathrm{Pr}\left( \left\Vert \overline{\psi}^{[1]}_{\acute{\otimes}} \right\Vert > \frac{\theta}{\left\Vert \frac{d X(t)}{dt} \right\Vert } \right) \nonumber \\
& \leq & \mathrm{Pr}\left( \left\Vert \overline{\psi}^{[1]}_{\acute{\otimes}} \right\Vert > \frac{\theta}{\Upsilon_X} \right),
\end{eqnarray}
where $\leq_1$ comes from Lemma~\ref{lma:norm estimation}. This lemma follows by Markov inequality.
$\hfill \Box$

Following lemma is about the expression of the $k$-th operator-valued function by MOIs. 
\begin{lemma}\label{lma:k-th derivative by k-MOI}
Let $A,B$ be two operators, and $\overline{\psi}_{\acute{\otimes}}$ be $k \in \mathbb{N}$ order differentiable function, then we have
\begin{eqnarray}\label{eq1:lma:k-th derivative by k-MOI}
\frac{d^k \overline{\psi}_{\acute{\otimes}}(A+tB)}{dt^k}\Bigg\vert_{t=0}&=&
k! T^{\overbrace{A+tB, \cdots, A+tB}^{k+1}}_{\overline{\psi}^{[k]}_{\acute{\otimes}}}\left(\underbrace{B,\cdots,B}_{k} \right).
\end{eqnarray}
\end{lemma}
\textbf{Proof:}
We will prove this lemma by induction. The base case is proved by Lemma~\ref{lma:1-th derivative by DOIs} for $k=1$. Suppose we have Eq.~\eqref{eq1:lma:k-th derivative by k-MOI} for the $k-1$-th derivative, then
\begin{eqnarray}\label{eq2:lma:k-th derivative by k-MOI}
\frac{d^k \overline{\psi}_{\acute{\otimes}}(A+tB)}{dt^k}\Bigg\vert_{t=0}&=&
\lim\limits_{t \rightarrow 0} \frac{(k-1)!}{t}\left( T^{\overbrace{A+tB, \cdots, A+tB}^{k}}_{\overline{\psi}^{[k-1]}_{\acute{\otimes}}}\left(\underbrace{B,\cdots,B}_{k-1} \right) - \right. \nonumber \\
&& \left. T^{\overbrace{A, \cdots, A}^{k}}_{\overline{\psi}^{[k-1]}_{\acute{\otimes}}}\left(\underbrace{B,\cdots,B}_{k-1} \right) \right)  \nonumber \\
&=& \lim\limits_{t \rightarrow 0} \frac{(k-1)!}{t}\sum\limits_{i=0}^{k-1}\left( T^{\overbrace{A+tB, \cdots, A+tB}^{k-i},\overbrace{A, \cdots, A}^{i}}_{\overline{\psi}^{[k-1]}_{\acute{\otimes}}}\left(\underbrace{B,\cdots,B}_{k-1} \right) - \right. \nonumber \\
&& \left. T^{\overbrace{A+tB, \cdots, A+tB}^{k-i-1},\overbrace{A, \cdots, A}^{i+1}}_{\overline{\psi}^{[k-1]}_{\acute{\otimes}}}\left(\underbrace{B,\cdots,B}_{k-1} \right) \right) \nonumber \\
&=_1& k!T^{\overbrace{A+tB, \cdots, A+tB}^{k-i},\overbrace{A, \cdots, A}^{i+1}}_{\overline{\psi}^{[k]}_{\acute{\otimes}}}\left(\underbrace{B,\cdots,B}_{k} \right),
\end{eqnarray}
where we apply perturbation formula, Lemma~\ref{lma:MOI perturbation}, in $=_1$.

From continuity property given by Lemma~\ref{lma:MOI continuity}, we have 
\begin{eqnarray}\label{eq3:lma:k-th derivative by k-MOI}
\lim\limits_{t \rightarrow 0} T^{\overbrace{A+tB, \cdots, A+tB}^{k-i},\overbrace{A, \cdots, A}^{i+1}}_{\overline{\psi}^{[k]}_{\acute{\otimes}}}\left(\underbrace{B,\cdots,B}_{k} \right) =  T^{\overbrace{A, \cdots, A}^{k+1}}_{\overline{\psi}^{[k]}_{\acute{\otimes}}}\left(\underbrace{B,\cdots,B}_{k} \right).
\end{eqnarray}

This lemma is proved from Eq.~\eqref{eq2:lma:k-th derivative by k-MOI} and Eq.~\eqref{eq3:lma:k-th derivative by k-MOI}.
$\hfill \Box$

From Lemma~\ref{lma:k-th derivative by k-MOI}, we are ready to have the following theorem about the tail bound for $\left\Vert \frac{d^k \overline{\psi}_{\acute{\otimes}}(A+tB)}{dt^k}\Bigg\vert_{t=0} \right\Vert$.
\begin{theorem}\label{thm:tail bound for the k der}
Let $A$ be random operators, $B$ be a determinstic operator, and $\overline{\psi}_{\acute{\otimes}}$ be $k \in \mathbb{N}$ order differentiable function, then we have 
\begin{eqnarray}\label{eq1:thm:tail bound for the k der}
\mathrm{Pr}\left( \left\Vert \frac{d^k \overline{\psi}_{\acute{\otimes}}(A+tB)}{dt^k}\Bigg\vert_{t=0} \right\Vert > \theta \right) & \leq & \frac{k! \left\Vert B \right\Vert^k   }{\theta} \mathbb{E}\left[ \left\Vert \overline{\psi}^{[k]}_{\acute{\otimes}} \right\Vert  \right].
\end{eqnarray}
\end{theorem}
\textbf{Proof:}
Because we have
\begin{eqnarray}
\mathrm{Pr}\left( \left\Vert \frac{d^k \overline{\psi}_{\acute{\otimes}}(A+tB)}{dt^k}\Bigg\vert_{t=0} \right\Vert > \theta \right) & =_1 & \mathrm{Pr}\left( \left\Vert k! T^{\overbrace{A+tB, \cdots, A+tB}^{k+1}}_{\overline{\psi}^{[k]}_{\acute{\otimes}}}\left(\underbrace{B,\cdots,B}_{k} \right) \right\Vert > \theta \right) \nonumber \\
&\leq_2 & \mathrm{Pr}\left(  k! \left\Vert \overline{\psi}^{[k]}_{\acute{\otimes}} \right\Vert \left\Vert B \right\Vert^k > \theta \right) \nonumber \\
&\leq_3& \frac{k! \left\Vert B \right\Vert^k   }{\theta} \mathbb{E}\left[ \left\Vert \overline{\psi}^{[k]}_{\acute{\otimes}} \right\Vert  \right],
\end{eqnarray}
where $=_1$ comes from Lemma~\ref{lma:k-th derivative by k-MOI}, $\leq_2$ comes from Lemma~\ref{lma:norm estimation}, and $\leq_3$ comes from Markov inequality again. 
$\hfill \Box$

\subsection{Tail Bound for Higher Random Operator Difference}\label{sec:Tail Bound for Higher Random Operator Difference}

In this section, we will extend the operator first order difference, i.e., $\overline{\psi}_{\acute{\otimes}}(A) - \overline{\psi}_{\acute{\otimes}}(A+B)$, to higher order difference. The higher opertor difference is defined as 
\begin{eqnarray}\label{eq:higher difference def}
\Delta^k_B \overline{\psi}_{\acute{\otimes}}(A)&=& \sum\limits_{i=0}^k (-1)^{k-i}{k \choose i}  \overline{\psi}_{\acute{\otimes}}(A+ i B),
\end{eqnarray} 
where $A,B$ are self-adjoint operators on Hilbert space. We will have 
the following theorem about the tail bound for the higher opertor difference given by Eq.~\eqref{eq:higher difference def}. 

\begin{theorem}\label{thm:tail bound for higher difference}
Let $A$ be a self-adjoint random operator, $B$ be a deterministic self-adjoint operator and the eigenvalue for the operator $A+iB$ is represented by $\lambda_{i+1}$ for $i=0,1,\cdots, k$. We assume that $\max(\lambda_{j+1} - \lambda_j) < \kappa$ for $j=1,2,\cdots, k$. Then, we have
\begin{eqnarray}\label{eq1:thm:tail bound for higher difference}
\mathrm{Pr}\left( \left\Vert \Delta^k_B \overline{\psi}_{\acute{\otimes}}(A) \right\Vert > \theta \right) &\leq& \frac{k \kappa \left\Vert B\right\Vert^k}{\theta} \mathbb{E}(\overline{\psi}_{\acute{\otimes}}).
\end{eqnarray}
\end{theorem}
\textbf{Proof:}
From the higher opertor difference provided by Eq.~\eqref{eq:higher difference def} and the divide difference definition, we have
\begin{eqnarray}\label{eq2:thm:tail bound for higher difference}
\Delta^k_B \overline{\psi}_{\acute{\otimes}}(A) &=& (-1)^{k}\Delta^{k-1}_B \overline{\psi}_{\acute{\otimes}}(A) + (-1)^{k-1}\Delta^{k-1}_B \overline{\psi}_{\acute{\otimes}}(A+B) \nonumber \\
&=& \sum\limits_{j=1}^{k}\left( \underbrace{\int \cdots \int}_{k+1} (\lambda_{j+1} - \lambda_{j})\overline{\psi}^{[k]}_{\acute{\otimes}}(\lambda_1, \cdots, \lambda_{k+1}) d P_{A}(\lambda_1)B d P_{A+B}(\lambda_2)B \right. \nonumber \\
&& \left.\cdots 
Bd P_{A+kB}(\lambda_{k+1}) \right),
\end{eqnarray}
where $P_A(\lambda_1), \cdots, d P_{A+kB}(\lambda_{k+1})$ are spectral measures. 

Then, we have
\begin{eqnarray}\label{eq3:thm:tail bound for higher difference}
\lefteqn{\mathrm{Pr}\left( \left\Vert \Delta^k_B \overline{\psi}_{\acute{\otimes}}(A) \right\Vert > \theta \right)} \nonumber \\
& =_1& \mathrm{Pr}\left( \left\Vert \sum\limits_{j=1}^{k}\left( \underbrace{\int \cdots \int}_{k+1} (\lambda_{j+1} - \lambda_{j})\overline{\psi}^{[k]}_{\acute{\otimes}}(\lambda_1, \cdots, \lambda_{k+1}) d P_{A}(\lambda_1)B\cdots 
B d P_{A+kB}(\lambda_{k+1}) \right) \right\Vert > \theta \right) \nonumber \\
&\leq_2& \mathrm{Pr}\left( \sum\limits_{j=1}^{k} \max(\lambda_{j+1} - \lambda_{j}) \left\Vert \overline{\psi}^{[k]}_{\acute{\otimes}} \right\Vert \left\Vert B\right\Vert^k > \theta \right) \nonumber \\
&\leq_3& \mathrm{Pr}\left( k \kappa \left\Vert \overline{\psi}^{[k]}_{\acute{\otimes}} \right\Vert \left\Vert B\right\Vert^k > \theta \right),
\end{eqnarray}
where $=_1$ comes from Eq.~\eqref{eq2:thm:tail bound for higher difference}, $=_2$ comes from Lemma~\ref{lma:norm estimation}, and $\leq_3$ comes from assumptions about $\max(\lambda_{j+1} - \lambda_j)$. The proof of this theorem is complete from Markov inequality. 
$\hfill \Box$

\subsection{Tail Bound for Taylor Remainder of Random Operator-valued Functions}\label{sec:Tail Bound for Taylor Remainder of Random Operator-valued Functions}

In this section, we will derive tail bounds for Taylor remainder of random operator-valued functions. 
We will consider two types of random operators: self-adjoint and unitary. Let us begin with \emph{multi-index} notation definitions. 

Multi-indices are represented by vector symbols as $\bm{a}$ or $\bm{b}$:
\begin{eqnarray}
\bm{a} = (a_1, a_2, \cdots, a_n),~\bm{b} = (b_1, b_2, \cdots, b_n),
\end{eqnarray}
where $a_i, b_i \in \{0,1,2,\cdots\}$. If $\bm{a}$ is a multi-index, we have the following notation defnitions:
\begin{eqnarray}\label{eq:multi-index def}
\left\vert \bm{a} \right\vert &=& a_1 + a_2 + \cdots +a_n, \nonumber \\
\bm{a} ! &=& a_1! a_2! \cdots a_n!, \nonumber \\
\bm{X}^{\bm{a}}&=& X_1^{a_1}X_2^{a_2}\cdots X_n^{a_n}~\mbox{where $\bm{x} = (X_1,X_2,\cdots, X_n)$}.
\end{eqnarray}

\subsubsection{Random Self-Adjoint Operators}\label{sec:Random Self-Adjoint Operators}

All operators discussed in this section are assumed self-adjoint.  Given a multi-variate operator-valued function $f(X_1, X_2, \cdots, X_n)$ and perturbed operators $(H_1, H_2, \cdots, H_n)$, we define the partial derivative 
for $\partial^{\bm{a}}f(\bm{X} + \bm{t} \circ \bm{H})$ as
\begin{eqnarray}\label{eq:partial deri def}
\partial^{\bm{a}}f(\bm{X} + \bm{t} \circ \bm{H}) = \frac{\partial^{\left\vert\bm{a}\right\vert}f(\bm{X} + \bm{t} \circ \bm{H})}{\partial t_1^{a_1}\partial t_2^{a_2}\cdots\partial t_n^{a_n}}\,
\end{eqnarray}
where $\bm{t}$ is the vector of real variables $(t_1, t_2, \cdots, t_n)$ and $\circ$ is the Hadamard product. Then, we define the Taylor remainder for $f(X_1, X_2, \cdots, X_n)$ with perturbed operators $(H_1, H_2, \cdots, H_n)$ as 
\begin{eqnarray}\label{eq:Taylor rem def}
R_{k,f,\bm{X}}(\bm{H}) &=& f(\bm{X} + \bm{H}) - \sum\limits_{\left\vert \bm{a} \right\vert < k} \frac{\partial^{\bm{a}}f(\bm{X} + \bm{t} \circ \bm{H})  }{\bm{a}!}\bigg\vert_{\bm{t}=\bm{0}}, 
\end{eqnarray}
where $\bm{t}$ is the vector of real variables $(t_1, t_2, \cdots, t_n)$ and $\circ$ is the Hadamard product. 

If the function $f(X_1, X_2, \cdots, X_n)$ can be decomposed as the summation of function constructible by the limit of integral projective tensor product $\acute{\otimes}$, i.e., 
\begin{eqnarray}\label{eq:func decompose form}
f(X_1, X_2, \cdots, X_n) &=& \overline{\phi}_{\acute{\otimes},1}(X_1) +  \overline{\phi}_{\acute{\otimes},2}(X_2) + \cdots +  \overline{\phi}_{\acute{\otimes},n}(X_n),
\end{eqnarray}
the expression given by Eq.~\eqref{eq:Taylor rem def} can be further expressed as
\begin{eqnarray}\label{eq:Taylor rem def decompose}
R_{k,f,\bm{X}}(\bm{H}) &=& \sum\limits_{j=1}^n \left( \overline{\phi}_{\acute{\otimes}}(X_j+H_j) - 
\sum\limits_{\ell=0}^{k-1}\frac{1}{\ell !} \frac{d^{\ell}}{d t^\ell_j}\overline{\phi}_{\acute{\otimes},j}(X_j + t_j H_j)\bigg\vert_{t_j=0}\right).
\end{eqnarray}

Following lemma is given to express $R_{k,f,\bm{X}}(\bm{H})$ by MOIs. 
\begin{lemma}\label{lma:sd remainder by MOIs}
Let the function $f: \mathbb{R}^n \rightarrow \mathbb{R}$ can be decomposed as $f(\bm{X})=\sum\limits_{j=1}^n\overline{\phi}_{\acute{\otimes},j}(X_j)$, we have 
\begin{eqnarray}
R_{k,f,\bm{X}}(\bm{H}) &=& \sum\limits_{j=1}^n T^{X_j + H_j,\cdots,X_j}_{\overline{\phi}^{[k]}_{\acute{\otimes},j}}(\underbrace{H_j,\cdots,H_j}_{k}).
\end{eqnarray}
\end{lemma}
\textbf{Proof:}
From Lemma~\ref{lma:k-th derivative by k-MOI}, we have 
\begin{eqnarray}\label{eq1:lma:sd remainder by MOIs}
R_{k,\overline{\phi}_{\acute{\otimes},j},X_j}(H_j) &=& T^{X_j + H_j,X_j}_{\overline{\phi}^{[1]}_{\acute{\otimes},j}}(H_j) - \sum\limits_{\ell=1}^{k-1}T^{X_j,\cdots,X_j}_{\overline{\phi}^{[k]}_{\acute{\otimes},j}}(\underbrace{H_j,\cdots,H_j}_{k}).
\end{eqnarray}
By applying Lemma~\ref{lma:MOI perturbation} to $T^{X_j + H_j,X_j}_{\overline{\phi}^{[1]}_{\acute{\otimes},j}}(H_j)$ at the  R.H.S. of Eq.~\eqref{eq1:lma:sd remainder by MOIs}, we have
\begin{eqnarray}\label{eq2:lma:sd remainder by MOIs}
R_{k,\overline{\phi}_{\acute{\otimes},j},X_j}(H_j) &=& T^{X_j + H_j,X_j,X_j}_{\overline{\phi}^{[2]}_{\acute{\otimes},j}}(H_j) - \sum\limits_{\ell=2}^{k-1}T^{X_j,\cdots,X_j}_{\overline{\phi}^{[k]}_{\acute{\otimes},j}}(\underbrace{H_j,\cdots,H_j}_{k}).
\end{eqnarray}
This lemma is proved by repeating this procedure $k-1$ times and the decomposition relation given by Eq.~\eqref{eq:Taylor rem def decompose}.
$\hfill \Box$

Following theorem is about the tail bound for $R_{k,f,\bm{X}}(\bm{H})$.
\begin{theorem}\label{thm:tail bound for sd remainder}
For $j=1,2,\cdots, n$, let $X_j$ be a self-adjoint random operator, $H_j$ be a deterministic self-adjoint operator. The function $f: \mathbb{R}^n \rightarrow \mathbb{R}$ is assumed to be decomposed as $f(\bm{X})=\sum\limits_{j=1}^n\overline{\phi}_{\acute{\otimes},j}(X_j)$. Then, we have
\begin{eqnarray}\label{eq1:thm:tail bound for sd remainder}
\mathrm{Pr}\left( \left\Vert R_{k,f,\bm{X}}(\bm{H}) \right\Vert > \theta \right) &\leq& \sum\limits_{j=1}^{n}\frac{n\left\Vert H_j\right\Vert^k}{\theta} \mathbb{E}(\overline{\psi}^{[k]}_{\acute{\otimes},j}).
\end{eqnarray}
\end{theorem}
\textbf{Proof:}
Because we have
\begin{eqnarray}
\mathrm{Pr}\left( \left\Vert R_{k,f,\bm{X}}(\bm{H})  \right\Vert > \theta \right) & =_1 & \mathrm{Pr}\left( \left\Vert \sum\limits_{j=1}^n T^{X_j + H_j,\cdots,X_j}_{\overline{\phi}^{[k]}_{\acute{\otimes},j}}(\underbrace{H_j,\cdots,H_j}_{k}) \right\Vert > \theta \right) \nonumber \\
&\leq_2 &  \mathrm{Pr}\left(\sum\limits_{j=1}^n  \left\Vert  T^{X_j + H_j,\cdots,X_j}_{\overline{\phi}^{[k]}_{\acute{\otimes},j}}(\underbrace{H_j,\cdots,H_j}_{k}) \right\Vert > \theta \right) \nonumber \\
&\leq&  \sum\limits_{j=1}^n \mathrm{Pr}\left( \left\Vert  T^{X_j + H_j,\cdots,X_j}_{\overline{\phi}^{[k]}_{\acute{\otimes},j}}(\underbrace{H_j,\cdots,H_j}_{k}) \right\Vert > \frac{\theta}{n} \right) \nonumber \\
&\leq_3& \sum\limits_{j=1}^n  \frac{n \left\Vert H_j \right\Vert^k   }{\theta} \mathbb{E}\left[ \left\Vert \overline{\psi}^{[k]}_{\acute{\otimes},j} \right\Vert  \right],
\end{eqnarray}
where $=_1$ comes from Lemma~\ref{lma:sd remainder by MOIs}, $\leq_2$ is due to norm triangle inequality, and $\leq_3$  comes from Lemma~\ref{lma:norm estimation} and Markov inequality again. 
$\hfill \Box$

\subsubsection{Random Unitary Operators}\label{sec:Random Unitary Operators}

Given a multi-variate operator-valued function $f(X_1, X_2, \cdots, X_n)$ and perturbed operators $(H_1, H_2, \cdots, H_n)$, which are assumed to be bounded self-adjoint operators, we use $e^{\iota \bm{t}\circ\bm{H}}$ to represent the vector \\
$(e^{\iota t_1 H_1}, e^{\iota t_2 H_2}, \cdots, e^{\iota t_n H_n})$, 
where $\iota \define \sqrt{-1}$. We define the partial derivative for $\partial^{\bm{a}}f(e^{\iota \bm{t}\circ\bm{H}}\circ \bm{X})$ as
\begin{eqnarray}\label{eq:partial deri def U}
\partial^{\bm{a}}f(e^{\iota \bm{t}\circ\bm{H}}\circ \bm{X}) = \frac{\partial^{\left\vert\bm{a}\right\vert}f(e^{\iota \bm{t}\circ\bm{H}}\circ \bm{X})}{\partial t_1^{a_1}\partial t_2^{a_2}\cdots\partial t_n^{a_n}}.
\end{eqnarray}

We first define the Taylor remainder for $f(X_1, X_2, \cdots, X_n)$ with perturbed operators $(H_1, H_2, \cdots, H_n)$ as 
\begin{eqnarray}\label{eq:Taylor rem U}
Q_{k,f,\bm{X}}(\bm{H}) &=& f(e^{\iota \bm{H}}\circ \bm{X}) - \sum\limits_{\left\vert \bm{a} \right\vert < k} \frac{\partial^{\bm{a}}f(e^{\iota \bm{t}\circ\bm{H}}\circ \bm{X})}{\bm{a}!}\bigg \vert_{\bm{t}=\bm{0}}.
\end{eqnarray}
If the function $f$ can also be decomposed as Eq.~\eqref{eq:func decompose form}, we can rewrite Eq.~\eqref{eq:Taylor rem U} as
\begin{eqnarray}\label{eq:Taylor rem U decompose}
Q_{k,f,\bm{X}}(\bm{H}) &=& \sum\limits_{j=1}^n \left( \overline{\phi}_{\acute{\otimes}}(e^{\iota  H_j} X_j) - 
\sum\limits_{\ell=0}^{k-1}\frac{1}{\ell !} \frac{d^{\ell}}{d t^\ell_j}\overline{\phi}_{\acute{\otimes},j}(e^{\iota t_j  H_j} X_j)\bigg\vert_{t_j=0}\right).
\end{eqnarray}

Following lemma is given to express $R_{k,f,\bm{X}}(\bm{H})$ by MOIs. 
\begin{lemma}\label{lma:unitary remainder by MOIs}
Let the function $f: \mathbb{R}^n \rightarrow \mathbb{R}$ can be decomposed as $f(\bm{X})=\sum\limits_{j=1}^n\overline{\phi}_{\acute{\otimes},j}(X_j)$, we have 
\begin{eqnarray}\label{eq1:lma:unitary remainder by MOIs}
\lefteqn{Q_{k,f,\bm{X}}(\bm{H})=} \nonumber \\
&&  \sum\limits_{j=1}^n \left( \sum\limits_{\ell=1}^k  
 \sum\limits_{\shortstack{ $i_1,\cdots,i_{\ell} \geq 1$ \\ $\left\vert \bm{i} \right\vert = k$}} T^{e^{\iota H_j} X_j, X_j, \cdots, X_j}_{\overline{\phi}^{[\ell]}_{\acute{\otimes},j}}
\left(\underbrace{\sum\limits_{m=i_1}^{\infty}\frac{(\iota H_j)^m}{m!} X_j, \frac{(\iota H_j)^{i_2}}{i_2 !} X_j, \cdots, \frac{(\iota H_j)^{i_\ell}}{i_{\ell} !} X_j}_{\ell} \right)\right).
\end{eqnarray}
\end{lemma}
\textbf{Proof:}
We will prove this lemma by induction. For basic case with $k=1$ and any $j$ between $1$ and $n$, we have
\begin{eqnarray}\label{eq2:lma:unitary remainder by MOIs}
Q_{1,\overline{\phi}_{\acute{\otimes},j},X_j}(H_j)&=&\overline{\phi}_{\acute{\otimes},j}(e^{\iota H_j} X_j) - \overline{\phi}_{\acute{\otimes},j}(X_j) \nonumber \\
&=& T^{e^{\iota H_j} X_j, X_j}_{\overline{\phi}^{[1]}_{\acute{\otimes},j}}(e^{\iota H_j} X_j - X_j) \nonumber \\
&=&  T^{e^{\iota H_j} X_j, X_j}_{\overline{\phi}^{[1]}_{\acute{\otimes},j}}\left(\sum\limits_{m=1}^{\infty}\frac{(\iota H_j)^{m}}{m!}X_j\right).
\end{eqnarray}
Therefore, by summing over $j$ with respect to Eq.~\eqref{eq2:lma:unitary remainder by MOIs}, this lemma is valid for the basic case.

Suppose that Eq.~\eqref{eq1:lma:unitary remainder by MOIs} is valid for $k=q$ and any $j$ between $1$ and $n$, we have
\begin{eqnarray}
\lefteqn{Q_{q+1,\overline{\phi}_{\acute{\otimes},j},X_j}(H_j)}
\nonumber \\
&=& \sum\limits_{\ell=1}^{q} 
\sum\limits_{\shortstack{ $i_1,\cdots,i_{\ell} \geq 1$ \\ $\left\vert \bm{i} \right\vert = q$}}
 T^{e^{\iota H_j} X_j, X_j, \cdots, X_j}_{\overline{\phi}^{[\ell]}_{\acute{\otimes},j}} \left(\underbrace{\sum\limits_{m=i_1}^{\infty}\frac{(\iota H_j)^m}{m!} X_j, \frac{(\iota H_j)^{i_2}}{i_2 !} X_j, \cdots, \frac{(\iota H_j)^{i_\ell}}{i_{\ell} !} X_j}_{\ell} \right) \nonumber \\
&&
- \frac{1}{q!} \frac{d^q}{d t^q_j}\overline{\phi}_{\acute{\otimes},j}(e^{\iota H_j} X_j)\bigg \vert_{t_j = 0} \nonumber 
\end{eqnarray}
\begin{eqnarray}\label{eq3:lma:unitary remainder by MOIs}
&=& \sum\limits_{\ell=1}^{q} 
\sum\limits_{\shortstack{ $i_1,\cdots,i_{\ell} \geq 1$ \\ $\left\vert \bm{i} \right\vert = q$}}
 T^{e^{\iota H_j} X_j, X_j, \cdots, X_j}_{\overline{\phi}^{[\ell]}_{\acute{\otimes},j}} \left(\underbrace{\sum\limits_{m=i_1}^{\infty}\frac{(\iota H_j)^m}{m!} X_j, \frac{(\iota H_j)^{i_2}}{i_2 !} X_j, \cdots, \frac{(\iota H_j)^{i_\ell}}{i_{\ell} !} X_j}_{\ell} \right) \nonumber \\
&&
-  \sum\limits_{\ell=1}^{q} 
\sum\limits_{\shortstack{ $i_1,\cdots,i_{\ell} \geq 1$ \\ $\left\vert \bm{i} \right\vert = q$}}
 T^{X_j, X_j, \cdots, X_j}_{\overline{\phi}^{[\ell]}_{\acute{\otimes},j}} \left(\underbrace{\frac{(\iota H_j)^{i_1}}{i_1 !} X_j, \frac{(\iota H_j)^{i_2}}{i_2 !} X_j, \cdots, \frac{(\iota H_j)^{i_\ell}}{i_{\ell} !} X_j}_{\ell} \right) \nonumber \\
&=& \sum\limits_{\ell=1}^{q} 
\sum\limits_{\shortstack{ $i_1,\cdots,i_{\ell} \geq 1$ \\ $\left\vert \bm{i} \right\vert = q$}}
\Bigg[T^{e^{\iota H_j} X_j, X_j, \cdots, X_j}_{\overline{\phi}^{[\ell]}_{\acute{\otimes},j}} \left(\underbrace{\sum\limits_{m=i_1}^{\infty}\frac{(\iota H_j)^m}{m!} X_j, \frac{(\iota H_j)^{i_2}}{i_2 !} X_j, \cdots, \frac{(\iota H_j)^{i_\ell}}{i_{\ell} !} X_j}_{\ell} \right) \nonumber \\
&&  -  T^{e^{\iota H_j}X_j, X_j, \cdots, X_j}_{\overline{\phi}^{[\ell]}_{\acute{\otimes},j}} \left(\underbrace{\frac{(\iota H_j)^{i_1}}{i_1 !} X_j, \frac{(\iota H_j)^{i_2}}{i_2 !} X_j, \cdots, \frac{(\iota H_j)^{i_\ell}}{i_{\ell} !} X_j}_{\ell} \right) \nonumber \\
&&
+
 T^{e^{\iota H_j}X_j, X_j, \cdots, X_j}_{\overline{\phi}^{[\ell]}_{\acute{\otimes},j}} \left(\underbrace{\frac{(\iota H_j)^{i_1}}{i_1 !} X_j, \frac{(\iota H_j)^{i_2}}{i_2 !} X_j, \cdots, \frac{(\iota H_j)^{i_\ell}}{i_{\ell} !} X_j}_{\ell} \right) \nonumber \\
&&-  T^{X_j, X_j, \cdots, X_j}_{\overline{\phi}^{[\ell]}_{\acute{\otimes},j}} \left(\underbrace{\frac{(\iota H_j)^{i_1}}{i_1 !} X_j, \frac{(\iota H_j)^{i_2}}{i_2 !} X_j, \cdots, \frac{(\iota H_j)^{i_\ell}}{i_{\ell} !} X_j}_{\ell} \right)
\Bigg] \nonumber \\
&=_1&  \overbrace{\sum\limits_{\ell=1}^{q}\sum\limits_{\shortstack{ $i_1,\cdots,i_{\ell} \geq 1$ \\ $\left\vert \bm{i} \right\vert = q$}}   T^{e^{\iota H_j} X_j, X_j, \cdots, X_j}_{\overline{\phi}^{[\ell]}_{\acute{\otimes},j}} \left(\underbrace{\sum\limits_{m=i_1 +1}^{\infty}\frac{(\iota H_j)^m}{m!} X_j, \frac{(\iota H_j)^{i_2}}{i_2 !} X_j, \cdots, \frac{(\iota H_j)^{i_\ell}}{i_{\ell} !} X_j}_{\ell} \right)}^{\mbox{Sum}_1} \nonumber \\
&& +  \overbrace{\sum\limits_{\ell=1}^{q}\sum\limits_{\shortstack{ $i_1,\cdots,i_{\ell} \geq 1$ \\ $\left\vert \bm{i} \right\vert = q$}}   T^{e^{\iota H_j} X_j, X_j, \cdots, X_j}_{\overline{\phi}^{[\ell+1]}_{\acute{\otimes},j}} \left(\underbrace{\sum\limits_{m=1}^{\infty}\frac{(\iota H_j)^m}{m!} X_j,\frac{(\iota H_j)^{i_1}}{i_1 !} X_j,  \frac{(\iota H_j)^{i_2}}{i_2 !} X_j, \cdots, \frac{(\iota H_j)^{i_\ell}}{i_{\ell} !} X_j}_{\ell+1} \right)}^{\mbox{Sum}_2} \nonumber \\
\end{eqnarray}
where $=_1$ comes from Lemma~\ref{lma:MOI perturbation}. 

By changing the summation index range as $\alpha_1 = i_1 + 1$ and $\alpha_p = i_p$ for $2 \leq p \leq \ell$, we can express the summation of $\mbox{Sum}_1$ as 
\begin{eqnarray}\label{eq4:lma:unitary remainder by MOIs}
\lefteqn{\mbox{Sum}_1=} \nonumber \\
&& \sum\limits_{\ell=1}^{q}\sum\limits_{\shortstack{ $\alpha_1 \geq 2,\alpha_2, \cdots, \alpha_{\ell} \geq 1$ \\ $\left\vert \bm{\alpha} \right\vert = q$}}   T^{e^{\iota H_j} X_j, X_j, \cdots, X_j}_{\overline{\phi}^{[\ell]}_{\acute{\otimes},j}} \left(\underbrace{\sum\limits_{m=\alpha_1}^{\infty}\frac{(\iota H_j)^m}{m!} X_j, \frac{(\iota H_j)^{\alpha_2}}{\alpha_2 !} X_j, \cdots, \frac{(\iota H_j)^{\alpha_\ell}}{i_{\ell} !} X_j}_{\ell} \right). 
\end{eqnarray}
Similarly, by changing the summation index range from $\ell=1,2,\cdots, q$ to $\ell=2,\cdots, q+1$, $\alpha_1 = 1$ and $\alpha_p = i_{p-1}$ for $2 \leq p \leq \ell$, we can express the summation of $\mbox{Sum}_2$ as 
\begin{eqnarray}\label{eq5:lma:unitary remainder by MOIs}
\lefteqn{\mbox{Sum}_2=} \nonumber \\
&& \sum\limits_{\ell=2}^{q+1}\sum\limits_{\shortstack{ $\alpha_1 =1, \alpha_2,\cdots,\alpha_{\ell} \geq 1$ \\ $\left\vert \bm{\alpha} \right\vert = q+1$}}  T^{e^{\iota H_j} X_j, X_j, \cdots, X_j}_{\overline{\phi}^{[\ell]}_{\acute{\otimes},j}} \left(\underbrace{\sum\limits_{m=\alpha_1}^{\infty}\frac{(\iota H_j)^m}{m!} X_j, \frac{(\iota H_j)^{\alpha_2}}{\alpha_2 !} X_j, \cdots, \frac{(\iota H_j)^{\alpha_\ell}}{\alpha_{\ell} !} X_j}_{\ell} \right). 
\end{eqnarray}

From Eqs.~\eqref{eq3:lma:unitary remainder by MOIs}~\eqref{eq4:lma:unitary remainder by MOIs} and~\eqref{eq5:lma:unitary remainder by MOIs}, we have
\begin{eqnarray}\label{eq6:lma:unitary remainder by MOIs}
\lefteqn{Q_{q+1,\overline{\phi}_{\acute{\otimes},j},X_j}(H_j)}
\nonumber \\
&=&  T^{e^{\iota H_j} X_j, X_j}_{\overline{\phi}^{[1]}_{\acute{\otimes},j}}\left(\sum\limits_{m=q+1}^{\infty}\frac{(\iota H_j)^{m}}{m!}X_j \right) \nonumber \\
&& +\sum\limits_{\ell=2}^{q+1}\left(\sum\limits_{\shortstack{ $\alpha_1 =1, \alpha_2,\cdots,\alpha_{\ell} \geq 1$ \\ $\left\vert \bm{\alpha} \right\vert = q+1$}} +
\sum\limits_{\shortstack{ $\alpha_1 \geq 2, \alpha_2,\cdots,\alpha_{\ell} \geq 1$ \\ $\left\vert \bm{\alpha} \right\vert = q+1$}} \right) \nonumber \\
&& T^{e^{\iota H_j} X_j, X_j, \cdots, X_j}_{\overline{\phi}^{[\ell]}_{\acute{\otimes},j}} \left(\underbrace{\sum\limits_{m=\alpha_1}^{\infty}\frac{(\iota H_j)^m}{m!} X_j, \frac{(\iota H_j)^{\alpha_2}}{\alpha_2 !} X_j, \cdots, \frac{(\iota H_j)^{\alpha_\ell}}{\alpha_{\ell} !} X_j}_{\ell} \right)  \nonumber \\
&&+ T^{e^{\iota H_j} X_j, X_j, \cdots, X_j}_{\overline{\phi}^{[q+1]}_{\acute{\otimes},j}} \left(\underbrace{\sum\limits_{m=1}^{\infty}\frac{(\iota H_j)^m}{m!} X_j, \iota H_j  X_j, \cdots, \iota H_j  X_j}_{q+1} \right) \nonumber \\
&=& \sum\limits_{\ell=1}^{q+1}
 \sum\limits_{\shortstack{ $i_1,\cdots,i_{\ell} \geq 1$ \\ $\left\vert \bm{i} \right\vert = q+1$}} T^{e^{\iota H_j} X_j, X_j, \cdots, X_j}_{\overline{\phi}^{[\ell]}_{\acute{\otimes},j}}
\left(\underbrace{\sum\limits_{m=i_1}^{\infty}\frac{(\iota H_j)^m}{m!} X_j, \frac{(\iota H_j)^{i_2}}{i_2 !} X_j, \cdots, \frac{(\iota H_j)^{i_\ell}}{i_{\ell} !} X_j}_{\ell} \right).
\end{eqnarray}
Then, this lemma is proved by induction and summing Eq.~\eqref{eq6:lma:unitary remainder by MOIs} with respect to the variable $j$.
$\hfill \Box$


We require another lemma to bound MOIs associated to $Q_{k,f,\bm{X}}(\bm{H})$. 
\begin{lemma}\label{lma:bound of MOI unitary}
Let $X_j$ be a random unitary operators, $H_j$ be an operator with bounded norm and $\overline{\phi}_{\acute{\otimes}, j}$ be a function constructible by the \emph{limit of integral projective tensor product}, then we have
\begin{eqnarray}\label{eq1:lma:bound of MOI unitary}
\left\Vert T^{e^{\iota H_j} X_j, X_j, \cdots, X_j}_{\overline{\phi}^{[\ell]}_{\acute{\otimes},j}}
\left(\underbrace{\sum\limits_{m=i_1}^{\infty}\frac{(\iota H_j)^m}{m!} X_j, \frac{(\iota H_j)^{i_2}}{i_2 !} X_j, \cdots, \frac{(\iota H_j)^{i_\ell}}{i_{\ell} !} X_j}_{\ell} \right) \right\Vert 
&\leq& \left\Vert \overline{\phi}^{[\ell]}_{\acute{\otimes},j} \right\Vert\rho(H_j, \bm{i}, \ell), 
\end{eqnarray}
where $\rho(H_j, \bm{i}, \ell)$ is defined as
\begin{eqnarray}\label{eq2:lma:bound of MOI unitary}
\rho(H_j, \bm{i}, \ell)&\define&  \left(\sum\limits_{m=i_1}^{\infty}\frac{\left\Vert H_j\right\Vert^m}{m!}\right)\prod\limits_{p=2}^{\ell}\frac{\left\Vert H_j\right\Vert^{i_p}}{i_p !}.
\end{eqnarray}
\end{lemma}
\textbf{Proof:}
Since we have
\begin{eqnarray}\label{eq3:lma:bound of MOI unitary}
\lefteqn{\left\Vert T^{e^{\iota H_j} X_j, X_j, \cdots, X_j}_{\overline{\phi}^{[\ell]}_{\acute{\otimes},j}}
\left(\underbrace{\sum\limits_{m=i_1}^{\infty}\frac{(\iota H_j)^m}{m!} X_j, \frac{(\iota H_j)^{i_2}}{i_2 !} X_j, \cdots, \frac{(\iota H_j)^{i_\ell}}{i_{\ell} !} X_j}_{\ell} \right) \right\Vert} \nonumber \\
&\leq_1& \left\Vert \overline{\phi}^{[\ell]}_{\acute{\otimes},j} \right\Vert \left\Vert \sum\limits_{m=i_1}^{\infty}\frac{(\iota H_j)^m}{m!} X_j \right\Vert \prod\limits_{p=2}^{\ell}\left\Vert \frac{(\iota H_j)^{i_p}}{i_p !} X_j\right\Vert  \nonumber \\
&\leq_2& \left\Vert \overline{\phi}^{[\ell]}_{\acute{\otimes},j} \right\Vert \left(\sum\limits_{m=i_1}^{\infty}\frac{\left\Vert H_j\right\Vert^m}{m!}\right)\prod\limits_{p=2}^{\ell}\frac{\left\Vert H_j\right\Vert^{i_p}}{i_p !} \nonumber \\
&=& \left\Vert \overline{\phi}^{[\ell]}_{\acute{\otimes},j} \right\Vert\rho(H_j,\bm{i},\ell), 
\end{eqnarray}
where $\leq_1$ comes from Lemma~\ref{lma:norm estimation} and $\leq_2$ comes from that the operator norm is submultiplicative and triangle inequality, this lemma is proved by Eq.~\eqref{eq2:lma:bound of MOI unitary}.
$\hfill \Box$

We are ready to present the following theorem about the tail bound for $Q_{k,f,\bm{X}}(\bm{H})$.
\begin{theorem}\label{thm:tail bound for unitary remainder}
For $j=1,2,\cdots, n$, let $X_j$ be a random unitary operator, $H_j$ be a deterministic operator with bounded norm. The function $f: \mathbb{R}^n \rightarrow \mathbb{R}$ is assumed to be decomposed as $f(\bm{X})=\sum\limits_{j=1}^n\overline{\phi}_{\acute{\otimes},j}(X_j)$. Then, we have
\begin{eqnarray}\label{eq1:thm:tail bound for unitary remainder}
\mathrm{Pr}\left( \left\Vert Q_{k,f,\bm{X}}(\bm{H}) \right\Vert > \theta \right) &\leq& \sum\limits_{j=1}^n  \sum\limits_{\ell=1}^k  \frac{kn \Theta(j,k,\ell)}{\theta}\mathbb{E}\left[ \left\Vert \overline{\psi}^{[\ell]}_{\acute{\otimes},j} \right\Vert  \right],
\end{eqnarray}
where $\Theta(j,k,\ell) \define \sum\limits_{\shortstack{ $i_1,\cdots,i_{\ell} \geq 1$ \\ $\left\vert \bm{i} \right\vert = k$}} \rho(H_j,\bm{i},\ell)$. Note that $\rho(H_j,\bm{i},\ell)$ is defined by Eq.~\eqref{eq2:lma:bound of MOI unitary}. 
\end{theorem}
\textbf{Proof:}
Because we have
\begin{eqnarray}
\lefteqn{\mathrm{Pr}\left( \left\Vert Q_{k,f,\bm{X}}(\bm{H})  \right\Vert > \theta \right)  =_1} \nonumber \\
&& \mathrm{Pr}\left( \left\Vert \sum\limits_{j=1}^n \left[ \sum\limits_{\ell=1}^k  
 \sum\limits_{\shortstack{ $i_1,\cdots,i_{\ell} \geq 1$ \\ $\left\vert \bm{i} \right\vert = k$}} \right. \right. \right. \nonumber \\
&&\left. \left. \left.  T^{e^{\iota H_j} X_j, X_j, \cdots, X_j}_{\overline{\phi}^{[\ell]}_{\acute{\otimes},j}}
\left(\underbrace{\sum\limits_{m=i_1}^{\infty}\frac{(\iota H_j)^m}{m!} X_j, \frac{(\iota H_j)^{i_2}}{i_2 !} X_j, \cdots, \frac{(\iota H_j)^{i_\ell}}{i_{\ell} !} X_j}_{\ell} \right)\right] \right\Vert > \theta \right) \nonumber \\
&\leq& \mathrm{Pr}\left(  \sum\limits_{j=1}^n  \sum\limits_{\ell=1}^k  
 \sum\limits_{\shortstack{ $i_1,\cdots,i_{\ell} \geq 1$ \\ $\left\vert \bm{i} \right\vert = k$}} \right. \nonumber \\
&&\left.    \left\Vert T^{e^{\iota H_j} X_j, X_j, \cdots, X_j}_{\overline{\phi}^{[\ell]}_{\acute{\otimes},j}}
\left(\underbrace{\sum\limits_{m=i_1}^{\infty}\frac{(\iota H_j)^m}{m!} X_j, \frac{(\iota H_j)^{i_2}}{i_2 !} X_j, \cdots, \frac{(\iota H_j)^{i_\ell}}{i_{\ell} !} X_j}_{\ell} \right) \right\Vert > \theta \right) \nonumber \\
&\leq_2& \mathrm{Pr}\left(  \sum\limits_{j=1}^n  \sum\limits_{\ell=1}^k  \left\Vert \overline{\psi}^{[\ell]}_{\acute{\otimes},j} \right\Vert
 \underbrace{\sum\limits_{\shortstack{ $i_1,\cdots,i_{\ell} \geq 1$ \\ $\left\vert \bm{i} \right\vert = k$}} \rho(H_j,\bm{i},\ell)}_{\define \Theta(j,k,\ell)}  > \theta \right) \nonumber \\
&\leq & \sum\limits_{j=1}^n  \sum\limits_{\ell=1}^k   \mathrm{Pr}\left( \left\Vert \overline{\psi}^{[\ell]}_{\acute{\otimes},j} \right\Vert  > \frac{\theta}{ kn \Theta(j,k,\ell) }\right) \nonumber \\
&\leq_3&  \sum\limits_{j=1}^n  \sum\limits_{\ell=1}^k  \frac{kn \Theta(j,k,\ell)}{\theta}\mathbb{E}\left[ \left\Vert \overline{\psi}^{[\ell]}_{\acute{\otimes},j} \right\Vert  \right],
\end{eqnarray}
where $=_1$ comes from Lemma~\ref{lma:unitary remainder by MOIs}, $\leq_2$ is due to Lemma~\ref{lma:bound of MOI unitary}, and $\leq_3$ Markov inequality again; this theorem is proved. 
$\hfill \Box$


\bibliographystyle{IEEETran}
\bibliography{RandomMOI_Bib}

\end{document}